\documentclass[aop,seceqn,citesort,MSNbibl,dvips]{arximspdf}
\usepackage[mathcal]{euscript}
\usepackage{mathrsfs}
\usepackage{graphics}
\doi{10.1214/09-AOP494}
\volume{38}
\issue{2}
\pubyear{2010}
\firstpage{498}
\lastpage{531}

\makeatletter

\DeclareMathAlphabet\mathcaligr{OMS}{cmsy}{m}{n}

\renewcommand{\mathscr}{\mathcal}

\newproclaim{definition}{Definition}[section]
\newtheorem{theorem}{Theorem}[section]
\newtheorem{lemma}{Lemma}[section]
\newtheorem{corollary}{Corollary}[section]
\newtheorem{proposition}{Proposition}[section]
\newtheorem{Step}{Step}

\makeatother

\begin{document}
\begin{frontmatter}

\title{A stochastic differential game for the inhomogeneous
$\infty$-Laplace equation}
\runtitle{SDG}

\begin{aug}
\author[A]{\fnms{Rami} \snm{Atar}\corref{}\ead[label=e1]{atar@ee.technion.ac.il}} and
\author[B]{\fnms{Amarjit} \snm{Budhiraja}\thanksref{t1}\ead[label=e2]{budhiraj@email.unc.edu}}
\runauthor{R. Atar and A. Budhiraja}
\affiliation{Technion---Israel Institute of Technology and University
of North Carolina}
\address[A]{Department of Electrical Engineering\\
Technion---Israel Institute of Technology\\
Haifa 32000\\
Israel\\
\printead{e1}} 
\address[B]{Department of Statistics\\
\quad and Operations Research\\
University of North Carolina\\
Chapel Hill, North Carolina 27599\\
USA\\
\printead{e2}}
\end{aug}

\thankstext{t1}{Supported in part by  Army Research Office
Grant W911NF-0-1-0080.}

\received{\smonth{8} \syear{2008}}
\revised{\smonth{1} \syear{2009}}

%
\begin{abstract}
Given a bounded $\mathcaligr{C}^2$ domain $G\subset{\mathbb{R}}^m$,
functions $g\in\mathcaligr{C}(\partial G,{\mathbb{R}})$ and $h\in\mathcaligr
{C}(\overline G,{\mathbb{R}}\setminus\{0\})$, let $u$ denote the unique
viscosity solution to the equation $-2\Delta_\infty u=h$ in $G$ with
boundary data $g$. We provide a representation for $u$ as the value of
a two-player zero-sum stochastic differential game.
\end{abstract}

%

\end{frontmatter}

\section{Introduction}\label{sec1}

\subsection{Infinity-Laplacian and games}

For an integer $m\ge2$,\vspace*{1pt} let a bounded $\mathcaligr{C}^2$ domain
$G\subset{\mathbb{R}}^m$, functions $g\in\mathcaligr{C}(\partial
G,{\mathbb{R}})$ and
$h\in\mathcaligr{C}(\overline G,{\mathbb{R}}\setminus\{0\})$ be given. We
study a two-player
zero-sum stochastic differential game (SDG), defined in terms of an
$m$-dimensional state process
that is driven by a one-dimensional Brownian motion, played until
the state exits the domain. The functions $g$ and $h$ serve as
terminal, and, respectively, running payoffs. The players' controls
enter in a diffusion coefficient and in an unbounded drift
coefficient of the state process. The dynamics are degenerate in
that it is possible for the players to completely switch off the
Brownian motion.
We show that the game has value, and characterize the value function
as the unique viscosity solution $u$ (uniqueness of solutions is
known from \cite{pssw}) of the equation
%
%
\begin{equation}\label{07}
\cases{-2\Delta_\infty u=h, &\quad in $G$,\cr
u=g, &\quad on $\partial G$.}
\end{equation}
Here, $\Delta_\infty$ is the infinity-Laplacian defined as $\Delta
_{\infty}f =
(Df)'(D^2f) (Df)/ |Df|^2$, provided $Df \neq0$, where for a
$\mathcaligr{C}^2$ function $f$ we denote by $Df$ the gradient and by $D^2f$
the Hessian matrix.
Our work is motivated by a representation for $u$ of Peres et al.
\cite{pssw} (established in fact in a far greater generality), as
the limit, as $\varepsilon\to0$, of the value function $V^\varepsilon
$ of a
discrete time random turn game, referred to as \textit{Tug-of-War}, in
which $\varepsilon$ is a parameter. The contribution of the current
work is
the identification of a game for which the value function is
precisely equal to $u$.

The infinity-Laplacian was first considered by Aronsson \cite{Aron}
in the study of absolutely minimal (AM) extensions of Lipschitz
functions. Given a Lipschitz function $u$ defined on the boundary
$\partial G$ of a domain $G$, a Lipschitz function $\widehat u$ extending
$u$ to $\overline G$ is called an AM extension of $u$ if, for every open
$U \subset G$, $\operatorname{Lip}_{\overline U}\widehat u = \operatorname
{Lip}_{\partial U}u$, where for
a real
function $f$ defined on $F \subset{\mathbb{R}}^m$,
$\operatorname{Lip}_F f = {\sup_{x,y \in F, x\neq y}} |f(x) - f(y)|/|x-y|$.
It was shown in \cite{Aron} that a Lipschitz function $\widehat u$ on
$\overline G$
that is $\mathcaligr{C}^2$
on $G$ is an AM extension of $\widehat u|_{\partial G}$ if and only if
$\widehat u$ is infinity-harmonic, namely satisfies $\Delta_{\infty
}\widehat u = 0$
in $G$.
This connection enables in some cases to prove uniqueness of AM
extensions via PDE tools.
However, due to the degeneracy of this elliptic equation, classical
PDE approach in general is not applicable.
Jensen \cite{Jen} showed that an appropriate framework
is through the theory of viscosity solutions,
by establishing existence and uniqueness of viscosity solutions to
the homogeneous version ($h=0$) of (\ref{07}), and showing that if
$g$ is Lipschitz then the solution is an AM extension of $g$. In
addition to the relation to AM extensions, the infinity-Laplacian
arises in a variety of other situations \cite{BEJ}. Some examples
include models for sand-pile evolution \cite{Aron2}, motion by mean
curvature and stochastic target problems \cite{KoSe,SoTo}.

We do not treat the homogenous equation for reasons mentioned later
in this section. The inhomogeneous equation may admit multiple
solutions when $h$ assumes both signs \cite{pssw}. Our assumption on
$h$ implies that either $h > 0 $ or $h < 0$. Uniqueness for the
case where these strict inequalities are replaced with weak inequalities
is unknown \cite{pssw}. Thus, the assumption we make on $h$ is the
minimal one under which
uniqueness is known to hold in general (except
the case $h=0$).

Let us describe the Tug-of-War game introduced in \cite{pssw}. Fix
$\varepsilon>0$.
Let a token be placed at $x\in G$, and set $X_0=x$. At the $k$th
step of the game ($k\ge1$), an independent toss of a fair coin
determines which player takes the turn. The selected player is
allowed to move the token from its current position $X_{k-1} \in G$
to a new position $X_k$ in $\overline G$, in such a way that $|X_k -
X_{k-1}| \le\varepsilon$ (\cite{pssw} requires $|X_k - X_{k-1}| <
\varepsilon$
but this is an equivalent formulation in the setting described
here).
The game ends at the first time $K$ when $X_K\in\partial G$. The
associated payoff is given by
%
%
\begin{equation}
\label{PF1051}
\mathbf{E} \Biggl[g(X_K) + \frac{\varepsilon^2}{4} \sum_{k=0}^{K-1}
h(X_k) \Biggr].
\end{equation}
Player I attempts to maximize the payoff and player II's goal is to
minimize it. It is shown in \cite{pssw} that
the value of the game, defined in a
standard way and denoted $V^\varepsilon(x)$, exists, that
$V^{\varepsilon}$
converges uniformly to a function $V$ referred to as the ``continuum
value function'' and that $V$ is the unique viscosity solution of
(\ref{07}) (these results are in fact also proved for the
homogeneous case, and in generality greater than the scope of the
current paper).
The question of associating a game directly with the continuum value
was posed and some basic technical challenges associated with it
were discussed in
\cite{pssw}.

Our approach to the question above is via a SDG formulation. To
motivate the form of the SDG, we start with the Tug-of-War game and
present some formal calculations (a precise definition of the SDG
will appear later). Let $\{\xi_k, k \in\mathbb{N}\}$ be a sequence
of i.i.d. random variables on some probability space $(\Omega,
\mathcaligr{F}, \mathbf{P})$ with $\mathbf{P}(\xi_k = 1) = \mathbf
{P}(\xi_k = -1)=1/2$,
interpreted as the sequence of coin tosses.
Let
$\{\mathcaligr{F}_k\}_{k\ge0}$ be a filtration of $\mathcaligr{F}$ to which
$\{\xi_k\}$ is adapted and such that $\{\xi_{k+1},\xi_{k+2},\ldots\}$
is independent of $\mathcaligr{F}_k$ for every $k\ge0$. Let $\{a_k\}$,
$\{b_k\}$ be $\{\mathcaligr{F}_k\}$-predictable
sequences of random variables with values in
$\overline{\mathbb{B}_{\varepsilon}(0)} = \{x \in{\mathbb{R}}^m\dvtx
|x| \le\varepsilon\}$. These
sequences correspond to control actions of players I and II;
that is, $a_k$ (resp., $b_k$) is the displacement exercised by
player I (resp., player II) if it wins the $k$th coin toss.
Associating the event $\{\xi_k = 1\}$ with player I winning the
$k$th toss, one can write the following representation for the
position of the token, starting from initial state $x$. For $j \in
\mathbb{N}$,
\[
X_j = x + \sum_{k=1}^j \biggl[
a_k\frac{1+\xi_k}{2}+b_k\frac{1-\xi_k}{2} \biggr]
=\sum_{k=1}^j\frac{a_k-b_k}{2}\xi_k+\sum_{k=1}^j\frac{a_k+b_k}{2}.
\]
We shall refer to $\{X_j\}$ as the ``state process.'' This
representation, in which turns are not taken at random
but both players select an action at each step, and the noise enters
in the dynamics, is more convenient for the development that
follows. Let $\varepsilon= 1/\sqrt{n}$ and rescale the control processes
by defining, for $t \ge0$, $A^n_t=\sqrt n  a_{[nt]}$, $B^n_t=\sqrt
n  b_{[nt]}$. Consider the continuous time state process $X^n_t =
X_{[nt]}$, and define
$\{W^n_t\}_{t\ge0}$ by setting $W^n_0=0$ and using the relation
\[
W^n_t = W^n_{(k-1)/n}+ \biggl(t-\frac{k-1}{n} \biggr)\sqrt{n}
\xi_k,\qquad
t\in \biggl(\frac{k-1}{n},  \frac{k}{n} \biggr],  k \in\mathbb{N}.
\]
Then we have
%
%
\begin{equation}\label{eq330} X^n_t = x + \frac{1}{2} \int_0^t
(A^n_s - B^n_s) \,dW^n_s + \frac{1}{2} \int_0^t \sqrt{n}(A^n_s +
B^n_s) \,ds.
\end{equation}
Note that $W^n$ converges weakly to a standard Brownian motion, and
since $|A^n_t| \vee|B^n_t| \le1$, the second term on the
right-hand side of (\ref{eq330}) forms a tight sequence. Thus, it is easy
to guess a substitute for it in the continuous game.
Interpretation of the asymptotics of the third term is more subtle,
and is a key element of the formulation.
One possible approach is to
replace the factor $\sqrt{n}$ by a large quantity that is
dynamically controlled by the two players. This point of view
motivates one to consider the identity (that we prove in Proposition \ref{prop2})
%
%
\begin{eqnarray}\label{44}
&&-2\Delta_\infty f = \sup_{|b|=1,  d \ge0 }   \inf_{|a|=1,  c\ge0}
\biggl\{-\frac{1}{2} (a-b)' (D^2f)  (a-b)\nonumber\\
&&\hspace*{148.4pt}{} - (c+d)(a+b)\cdot Df
\biggr\},\\
\eqntext{f \in\mathcaligr{C}^2, Df \neq0,}
\end{eqnarray}
for the following reason.
Let $\mathcaligr{H} = \mathcaligr{S}^{m-1} \times[0, \infty)$
where $\mathcaligr{S}^{m-1}$ is the unit sphere in ${\mathbb{R}}^m$. The
expression in curly brackets is equal to $\mathcaligr{L}^{a,b,c,d}f(x)$,
where for $(a,c),(b,d)\in\mathcaligr{H}$, $\mathcaligr{L}^{a,b,c,d}$ is
the controlled
generator associated with the process
%
%
\begin{equation}\label{star1050}
X_t = x+\int_0^t(A_s-B_s)\,dW_s+\int_0^t(C_s+D_s)(A_s+B_s)\,ds,
\qquad
t\in[0,\infty),\hspace*{-33pt}
\end{equation}
and $(A, C)$ and $(B, D)$ are control processes taking values in
$\mathcaligr{H}$.
Since $\Delta_\infty$ is related to (\ref{eq330}) via the
Tug-of-War, and
$\mathcaligr{L}^{a,b,c,d}$ to (\ref{star1050}), identity (\ref{44}) suggests
to regard (\ref{star1050}) as a formal limit of (\ref{eq330}).
Consequently the SDG will have (\ref{star1050}) as a state process,
where the controls $(A,C)$ and $(B,D)$ are chosen by the two players.
Finally, the payoff functional, as a formal limit of (\ref{PF1051}), and
accounting for the extra factor of $1/2$ in (\ref{eq330}), will be
given by
$\mathbf{E}[\int_0^{\tau} h(X_s) \,ds + g(X_{\tau})]$, where
$\tau=\inf\{t\dvtx X_t \notin G \}$
(with an appropriate convention regarding $\tau= \infty$).

A precise formulation of this game is given in Section \ref{sec2},
along with a statement of the main result.
Section \ref{sec1.3} discusses the technique and some open problems.

Throughout, we will denote by $\mathscr{S}(m)$ the space of symmetric $m
\times m$
matrices, and by $I_m\in\mathscr{S}(m)$ the identity matrix.
A function $\vartheta\dvtx[0,\infty)\to[0,\infty)$ will be said to be a
\textit{modulus} if it is continuous, nondecreasing, and satisfies
$\vartheta(0)=0$.

\subsection{SDG formulation and main result}\label{sec2}

Recall that $G$ is a bounded $\mathcaligr{C}^2$ domain in ${\mathbb
{R}}^m$, and that
$g\dvtx\partial G\to{\mathbb{R}}$ and $h\dvtx\overline G\to{\mathbb
{R}}\setminus\{0\}$ are given continuous
functions. In particular we have that either $h > 0$ or $h < 0$.
Since the two cases are similar, we will only consider $h>0$, and
use the notation $\underline h:=\inf_{\overline G} h>0$. Let
$(\Omega,\mathcaligr{F},\{\mathcaligr{F}_t\},\mathbf{P})$ be a complete
filtered probability
space with right-continuous filtration, supporting an
$(m+1)$-dimensional $\{\mathcaligr{F}_t\}$-Brownian motion $\overline
W=(W,\widetilde
W)$, where $W$ and $\widetilde W$ are one- and $m$-dimensional Brownian
motions, respectively. Let $\mathbf{E}$ denote expectation with
respect to
$\mathbf{P}$. Let $X_t$ be a process taking values in ${\mathbb
{R}}^m$, given by
%
%
\begin{equation}\label{01}
X_t = x+\int_0^t(A_s-B_s)\,dW_s+\int_0^t(C_s+D_s)(A_s+B_s)\,ds,\qquad
t\in[0,\infty),\hspace*{-33pt}
\end{equation}
where $x\in\overline G$, $A_t$ and $B_t$ take values in the unit sphere
$\mathcaligr{S}^{m-1}\subset{\mathbb{R}}^m$, and $C_t$ and $D_t$ take
values in
$[0,\infty)$. Denote
%
%
\begin{equation}\label{23}
Y^0=(A,C),\qquad  Z^0=(B,D).
\end{equation}
The processes $Y^0$ and $Z^0$ take values in
$\mathcaligr{H}=\mathcaligr{S}^{m-1}\times[0,\infty)$. These processes
will correspond to
control actions of the maximizing and minimizing player,
respectively. We remark that,\vspace*{1pt} although $\widetilde W$ does not appear
explicitly in the dynamics~(\ref{01}), the control processes $Y^0,
Z^0$ will be required to be $\{\mathcaligr{F}_t\}$-adapted, and thus may
depend on
it. In Section \ref{sec1.3}, we comment on the need for including
this auxiliary Brownian motion in our formulation. Let
\[
\tau=\inf\{t\dvtx X_t\in\partial G\}.
\]
Throughout, we will follow the convention that the infimum over an
empty set is~$\infty$. We write
%
%
\begin{equation}
\label{30}
X(x,Y^0,Z^0) \qquad [\mbox{resp., } \tau(x,Y^0,Z^0)]
\end{equation}
for the process $X$ (resp., the random time $\tau$) when it is
important to specify the explicit dependence on $(x,Y^0,Z^0)$. If
$\tau<\infty$ a.s., then the payoff $J(x, Y^0, Z^0)$ is well defined
with values in $(-\infty,\infty]$, where
%
%
\begin{equation}\label{02}
J(x,Y^0,Z^0)=\mathbf{E}\biggl[\int_0^\tau h(X_s)\,ds+g(X_\tau)\biggr]
\end{equation}
and $X$ is given by (\ref{01}). When $\mathbf{P}(\tau
(x,Y^0,Z^0)=\infty)>0$,
we set $J(x,Y^0,Z^0)=\infty$, in agreement with the expectation of the
first term in (\ref{02}).

We turn to the precise definition of the SDG.
For a process $H^0=(A,C)$ taking values in $\mathcaligr{H}$, we let
$S(H^0)=\operatorname{ess}\sup \sup_{t\in[0,\infty)}C_t$. In the
formulation below,
each player initially declares a bound $S$, and then plays so as to
keep $S(H^0)\le S$.
\begin{definition}\label{def1}
(i) A pair $H=(\{H^0_t\},S)$, where $S\in{\mathbb{N}}$ and $\{H^0_t\}
$ is a
process taking values in $\mathcaligr{H}$, is said to be an admissible
control if $\{H^0_t\}$ is $\{\mathcaligr{F}_t\}$-progressively measurable,
and
$S(H^0)\le S$.
The set of all admissible controls is denoted by
$M$.
For $H=(\{H^0_t\},S)\in M$, denote $\mathbf{S}(H)=S$.

(ii) A mapping $\varrho\dvtx M\to M$ is said to be a strategy if, for
every $t$,
\[
\mathbf{P}(H^0_s=\widetilde H^0_s \mbox{ for a.e. } s\in[0,t])=1
\quad\mbox{and}\quad
S=\widetilde S
\]
implies
\[
\mathbf{P}(I^0_s=\widetilde I^0_s \mbox{ for a.e. } s\in[0,t])=1
\quad\mbox{and}\quad
T=\widetilde T,
\]
where $(I^0,T)=\varrho[(H^0,S)]$ and $(\widetilde I^0,\widetilde
T)=\varrho[(\widetilde
H^0,\widetilde
S)]$.
The set of all strategies is denoted by $\widetilde\Gamma$. For
$\varrho
\in\widetilde\Gamma$, let $\mathbf{S}(\varrho) = \sup_{H \in
M}\mathbf{S}(\varrho[H])$.
Let
\[
\Gamma= \{\varrho\in\widetilde\Gamma\dvtx \mathbf{S}(\varrho)
<\infty\}.
\]
\end{definition}

We will use the symbols $Y$ and $\alpha$ for generic control and
strategy\break for the maximizing player, and $Z$ and $\beta$ for the
minimizing player. If\break
$Y = (Y^0,K), Z =(Z^0,L)\in M$, we sometimes write
$J(x,Y,Z) = J(x,(Y^0,\break K),(Z^0,L))$ for $J(x$,$Y^0,Z^0)$. Similar
conventions will be used for $X(x, Y, Z)$ and $\tau(x, Y, Z)$. Let
\begin{eqnarray*}
J^x(Y,\beta) &=& J(x,Y,\beta[Y]), \qquad x\in\overline G, Y\in M, \beta\in
\Gamma,
\\
J^x(\alpha,Z) &=& J(x,\alpha[Z],Z),\qquad  x\in\overline G, \alpha\in\Gamma
, Z\in M.
\end{eqnarray*}
Define analogously $X^x(Y,\beta)$, $X^x(\alpha,Z)$, $\tau^x(Y,\beta)$
and $\tau^x(\alpha,Z)$ via (\ref{30}). Define the lower value of the
SDG by
%
%
\begin{equation}\label{05}
V(x)=\inf_{\beta\in\Gamma}\sup_{Y\in M}J^x(Y,\beta)
\end{equation}
and the upper value by
%
%
\begin{equation}\label{06}
U(x)=\sup_{\alpha\in\Gamma}\inf_{Z\in M}J^x(\alpha,Z).
\end{equation}
The game is said to have a value if $U=V$.

Recall that the infinity-Laplacian is defined by $\Delta_\infty
f=p'\Sigma
p/|p|^2$, where $f$ is a $\mathcaligr{C}^2$ function, $p=Df$ and $\Sigma=D^2f$,
provided that $p\ne0$. Thus, $\Delta_\infty f$ is equal to the second
derivative in the direction of the gradient. In the special case
where $D^2f(x)$ is of the form $\lambda I_m$ for some real $\lambda$,
it is
therefore natural to define $\Delta_\infty f(x)=\lambda$ even if $Df(x)=0$
\cite{pssw}. This will be reflected in the definition of viscosity
solutions of (\ref{07}), that we state below.
Let
\begin{eqnarray*}
\mathcaligr{D}_0 &=& \{(0,\lambda I_m)\in{\mathbb{R}}^m\times\mathscr
{S}(m)\dvtx \lambda\in{\mathbb{R}}\},
\\
\mathcaligr{D}_1 &=& ({\mathbb{R}}^m\setminus\{0\})\times\mathscr{S}(m),
\\
\mathcaligr{D} &=& \mathcaligr{D}_0\cup\mathcaligr{D}_1
\end{eqnarray*}
and
\[
\Lambda(p,\Sigma)= \cases{-2\lambda, &\quad $(p,\Sigma)=(0,\lambda
I_m)\in
\mathcaligr{D}_0$, \vspace*{2pt}\cr
-2\dfrac{p'\Sigma p}{|p|^2}, &\quad $(p,\Sigma)\in\mathcaligr{D}_1$.}
\]
\begin{definition}\label{def3}
A continuous function $u\dvtx \overline G\to{\mathbb{R}}$ is said to be a viscosity
supersolution (resp., subsolution) of (\ref{07}), if:

\begin{longlist}
\item for every $x \in G$ and $\varphi\in\mathcaligr{C}^2(G)$ for which
$(p,\Sigma):=(D\varphi(x),D^2\varphi(x))\in\mathcaligr{D}$, and
$u - \varphi$ has a global minimum [maximum] on $G$ at $x$, one has
%
%
\begin{equation}\label{24}
\Lambda(p,\Sigma)-h(x)\ge0 \qquad [\le0];
\end{equation}
and
\item  $u=g$ on $\partial G$.

A viscosity solution is a function which is both a super- and a
subsolution.
\end{longlist}
\end{definition}

The result below has been established in \cite{pssw}.
\begin{theorem}\label{thpssw}
There exists a unique viscosity solution to (\ref{07}).
\end{theorem}

The following is our main result.
\begin{theorem}\label{th1}
The functions $U$ and $V$ are both viscosity solutions to
(\ref{07}). Consequently, the SDG has a value.
\end{theorem}

In what follows, we use the terms subsolution, supersolution and
solution as shorthand for viscosity subsolution, etc.

\subsection{Discussion}
\label{sec1.3}

We describe here our approach to proving the main result, and
mention some obstacles in extending it.

A common approach to showing solvability of Bellman--Isaacs (BI)
equations [(\ref{07}) can be viewed as such an equation due to
(\ref{44})] by the associated value function, is by proving that the
value function satisfies a dynamic programming principle (DPP).
Roughly speaking, this is an equation expressing the fact that,
rather than attempting to maximize their profit by considering
directly the payoff functional, the players may consider the payoff
incurred up to a time $t$ plus the value function evaluated at the
position $X_t$ that the state reaches at that time. Although in a
single player setting (i.e., in pure control problems) DPP are well
understood, game theoretic settings as in this paper are
significantly harder. In particular, as we shall shortly point out,
there are some basic open problems related to such DPP. In a setting
with a finite time horizon, Fleming and Souganidis \cite{FS}
established a DPP
based on careful discretization and approximation arguments. We have
been unable to carry out a similar proof in the current setting,
which includes a payoff given in terms of an exit time, degenerate
diffusion and unbounded controls.

Swiech \cite{swi} has developed an alternative approach to the above
problem that relies on existence of solutions. Instead of
establishing a DPP for the
value function, the idea of \cite{swi} is to show that
any \textit{solution} must satisfy a DPP. To see what is meant by such
a DPP and how it is used, consider the equation, $-2\Delta_\infty
u+\lambda
u=h$ in $G$, $u=g$ on $\partial G$, where $\lambda\ge0$ is a constant,
associated with the payoff in (\ref{02}) modified by a discount
factor. Assume that one can show that whenever $u$ and $v$ are sub-
and supersolutions, respectively, then
%
%
\begin{eqnarray}\label{45}
u(x)&\le&\sup_{\alpha\in\Gamma}\inf_{Z\in M}\mathbf{E} \biggl[\int
_0^\sigma
e^{-\lambda s}h(X_s)\,ds+e^{-\lambda\sigma}u(X_\sigma) \biggr],
\\
\label{46}
v(x)&\ge&\sup_{\alpha\in\Gamma}\inf_{Z\in M}\mathbf{E} \biggl[\int
_0^\sigma
e^{-\lambda s}h(X_s)\,ds+e^{-\lambda\sigma}v(X_\sigma) \biggr],
\end{eqnarray}
for $X=X[x,\alpha[Z],Z]$, $\tau=\tau[x,\alpha[Z],Z]$ and $\sigma
=\sigma(t)=\tau
\wedge
t$. Sending $t \to\infty$ in the above equations, one would
formally obtain
%
%
\begin{equation}\label{ins1155}u(x) \le
\sup_{\alpha\in\Gamma}\inf_{Z\in M}\mathbf{E} \biggl[\int_0^{\tau}
e^{-\lambda s}h(X_s)\,ds+e^{-\lambda\tau}g(X_{\tau}) \biggr] \le v(x),
\end{equation}
in particular yielding that if $u=v$ is a solution to the equation
then it must equal the upper value function. This would establish
unique solvability of the equation by the upper value function,
provided there exists a solution. In
the case $\lambda> 0$, justifying the above formal limit is
straightforward (see \cite{swi})
but the case $\lambda= 0$, as in our setting, requires a more careful
argument.
Our proofs exploit the uniform positivity of $h$ due to which the
minimizing player will not allow $\tau$ to be too large. This leads
to uniform estimates on the decay of $\mathbf{P}(\tau> t)$ as $t\to
\infty$,
from which an inequality as in (\ref{ins1155}) follows readily. This
discussion also explains why we are unable to treat the case $h=0$.

Establishing DPP as in (\ref{45}), (\ref{46}) is thus a key
ingredient in this approach. For a class of BI equations, defined
on all of ${\mathbb{R}}^m$, for which the associated game has a
bounded action
set and a fixed, finite time horizon, such a DPP was proved in
Swiech \cite{swi}. In the current paper, although we do not
establish (\ref{45}), (\ref{46}) in the above form, we derive
similar inequalities (for $\lambda= 0$) for a related bounded action
game, defined on $G$. The characterization of the value function
for the original unbounded action game is then treated by taking
suitable limits.

Both \cite{FS} and \cite{swi} require some assumptions on the
sample space and underlying filtration. In \cite{FS}, the underlying
filtration is the one generated by the driving Brownian motion. The
approach taken in \cite{swi}, which the current paper follows,
allows for a general filtration as long as it is rich enough to
support an $m$-dimensional Brownian motion, independent of the
Brownian motion driving the state process [for example, it could be
the filtration generated by an $(m+1)$-dimensional Brownian motion].
The reason for imposing this requirement in \cite{swi} is that
inequalities similar to (\ref{45}) and (\ref{46}) are proved by
first establishing them for a game associated with a nondegenerate
elliptic equation, and then taking a vanishing viscosity limit. This
technical issue is the reason for including the auxiliary process
$\widetilde W$ in our formulation as well. As pointed out in \cite{swi},
the question of validity of the DPP and the characterization of the
value as the unique solution to the PDE, under an arbitrary
filtration, remains a basic open problem on SDGs.

 The unboundedness of the action space, on one hand, and the combination of
degeneracy of the dynamics and an exit time criterion on the other hand, make it hard
to adapt the results of \cite{swi} to our setting. In order to
overcome the first difficulty, we approximate the original SDG by a
sequence of games with bounded action spaces, that are more readily
analyzed.
For the bounded action game, existence of solutions to the upper and
lower BI equations follow from \cite{CKLS}. We show that the
solutions to these equations satisfy a DPP similar to
(\ref{45}) and (\ref{46}) (Proposition \ref{prop1}). As discussed
above, existence of solutions along with the DPP yields the
characterization of these solutions as the corresponding value
functions. Next, as we show in Lemma \ref{lem5}, the upper and lower
value functions for the bounded action games approach the
corresponding value functions of the original game, pointwise, as
the bounds approach~$\infty$. Moreover, in Lemma \ref{lem4}, we
show that any uniform subsequential limit, as the bounds
approach $\infty$, of
solutions to the BI equation for bounded action games is a viscosity
solution of (\ref{07}). The last piece in the proof of the main
result is then showing existence of uniform (subsequential) limits.
This is established in Theorem \ref{th3} by proving equicontinuity,
in the parameters governing the bounds,
of the value functions for bounded action games. The proof of
equicontinuity is the most technical part of this paper and the main
place where the $\mathcaligr{C}^2$ assumption on the domain is used.
This is also the place where the possibility of degenerate dynamics
close to the exit time needs to be carefully analyzed.

The rest of this paper is organized as follows. In Section
\ref{sec3}, we prove Theorem~\ref{th1} based on results on BI
equations for bounded action SDG. These results are established in
Sections \ref{sec4} (equicontinuity of the value functions) and
\ref{sec5} (relating the value function to the PDE). Finally, it is
natural to ask whether the state process, obtained under
$\delta$-optimal play by both players, converges in law as $\delta$
tends to zero. Section \ref{sec6} describes a recently obtained
result \cite{AtBu2} that addresses this issue.

\section{Relation to Bellman--Isaacs equation}\label{sec3}

In this section, we prove Theorem \ref{th1} by relating the value
functions $U$ and $V$ to value functions of SDG with bounded action
sets, and similarly, the solution to (\ref{07}) to that of the
corresponding Bellman--Isaacs equations.

Let $p\in{\mathbb{R}}^m$, $p\ne0$ and $S\in\mathscr{S}(m)$ be
given, and, for
$n\in{\mathbb{N}}$, fix $p_n\in{\mathbb{R}}^m$, $p_n\ne0$ and
$S_n\in\mathscr{S}(m)$, such
that $p_n\to p$, $S_n\to S$. Denote $\overline p=p/|p|$ and $\overline
p_n=p_n/|p_n|$. Let $\{k_n\}$ and $\{l_n\}$ be positive, increasing
sequences such that $k_n\to\infty$, $l_n\to\infty$.

Denote
%
%
\begin{equation}\label{09}
\Phi(a,b,c,d;p,S)= -\tfrac12 (a-b)'S(a-b)-(c+d)(a+b)\cdot p,
\end{equation}
and let
%
%
\begin{eqnarray}\label{18}
\Lambda_{kl}^+(p,S)&=&
\max_{|b|=1, 0\le d\le l}  \min_{|a|=1, 0\le c\le k}
\Phi(a,b,c,d;p,S),
\\
\label{10}
\Lambda_{kl}^-(p,S)&=&\min_{|a|=1, 0\le c\le k}  \max_{|b|=1, 0\le
d\le
l} \Phi(a,b,c,d;p,S).
\end{eqnarray}
Set
\[
\Lambda_n^+(p,S)=\Lambda^+_{k_nl_n}(p,S),\qquad
\Lambda_n^-(p,S)=\Lambda^-_{k_nl_n}(p,S).
\]
\begin{lemma}\label{lem1}
One has $\Lambda^+_n(p_n,S_n)\to\Lambda(p,S)$, and
$\Lambda^-_n(p_n,S_n)\to\Lambda(p,S)$, as $n\to\infty$.
\end{lemma}
\begin{pf}
We prove only the statement regarding $\Lambda^-_n$, since the
other statement can be proved analogously. We omit the superscript
``$-$'' from the notation. Denote
$\Phi_n(a,b,c,d) = \Phi(a,b,c,d;p_n,S_n)$. Let
\[
\overline\Lambda_n(a,c)=\max_{|b|=1, 0\le d\le l_n}\Phi_n(a,b,c,d).
\]
Let $(a^*_n,c^*_n)$ be such that
$\Lambda_n^*:=\Lambda_n(p_n,S_n)=\overline\Lambda_n(a^*_n,c^*_n)$.
Note that $\Lambda_n^*\le\overline\Lambda_n(\overline p_n,0)$, which is
bounded from
above as $n\to\infty$, since $(b + \overline p_n)\cdot\overline p_n \ge0$ for
all $b \in\mathcaligr{S}^{m-1}$, $n \ge1$.
On the other hand, if for some fixed $\varepsilon>0$,
$a^*_n\cdot p_n<|p_n|-\varepsilon$ holds for infinitely many $n$, then
$\limsup\overline\Lambda_n(a^*_n,c_n)=\infty$ for any choice of $c_n$
contradicting the statement that $\Lambda_n^*$ is bounded from above.
This shows, for every $\varepsilon> 0$,
\[
|p_n|-\varepsilon\le a^*_n\cdot p_n\le|p_n|
\]
for all large $n$. In particular, $a^*_n\to\overline p$. Next note that
\[
\Lambda_n^*=\overline\Lambda_n(a^*_n,c^*_n)\ge\Phi_n(-\overline p_n,a^*_n,l_n,c^*_n)
\ge-\tfrac12 (\overline p_n+a^*_n)'S_n(\overline p_n+a^*_n)
\]
hence,
\[
\liminf\Lambda_n^*\ge-2\overline p'S\overline p=\Lambda(p,S).
\]
Also, with $(\widetilde b_n,\widetilde d_n)\in\arg\max_{(b,d)}\Phi
_n(b,\overline
p_n,d,k_n)$,
%
%
\begin{eqnarray}\label{star1106}
\Lambda_n^* &=& \overline\Lambda_n(a^*_n,c^*_n)\le\overline\Lambda_n(\overline
p_n,k_n) \nonumber\\
&=&-\tfrac12
(\widetilde b_n-\overline p_n)'S_n(\widetilde b_n-\overline p_n) - (\widetilde
d_n + k_n)(\widetilde b_n
+ \overline p_n) \cdot p_n
\\
&\le& -\tfrac12 (\widetilde b_n-\overline p_n)'S_n(\widetilde b_n-\overline
p_n).\nonumber
\end{eqnarray}
If $\widetilde b_n\to-\overline p$ does not hold, then $\liminf\Lambda
_n^*=-\infty$
by the first line of (\ref{star1106}) which
contradicts the previous display. This shows $\widetilde
b_n\to-\overline p$. Hence, from the second line of (\ref{star1106})
\[
\limsup\Lambda_n^*\le-2\overline p'S\overline p=\Lambda(p,S).
\]
\upqed\end{pf}

We now consider two formulations of SDG with bounded controls, the
first being based on Definition \ref{def1} whereas the second is more
standard.
For $k,l\in{\mathbb{N}}$,
let
\begin{eqnarray*}
M_k &=& \{Y\in M\dvtx\mathbf{S}(Y)\le k\},
\\
\Gamma_l &=& \{\beta\in\Gamma\dvtx\mathbf{S}(\beta)\le l\}.
\end{eqnarray*}
Define accordingly the lower value
%
%
\begin{equation}\label{20}
V_{kl}(x)=\inf_{\beta\in\Gamma_l}\sup_{Y\in M_k}J^x(Y,\beta),
\end{equation}
and the upper value
%
%
\begin{equation}\label{21}
U_{kl}(x)=\sup_{\alpha\in\Gamma_k}\inf_{Z\in M_l}J^x(\alpha,Z).
\end{equation}
\begin{definition}\label{def2}
(i) A process $\{H_t\}$ taking values in $\mathcaligr{H}$ is said to be a
simple admissible control if it is $\{\mathcaligr{F}_t\}$-progressively
measurable.
We denote by $M^0$ the set of all simple admissible controls, and
let $M^0_k=\{H\in M^0\dvtx S(H)\le k\}$.

(ii) Given $k,l\in{\mathbb{N}}$, we say that a mapping $\varrho
\dvtx M^0_k\to M^0_l$
is a simple strategy, and write $\varrho\in\Gamma^0_{kl}$ if, for every
$t$,
\[
\mathbf{P}(H_s=\widetilde H_s \mbox{ for a.e. } s\in[0,t])=1
\]
implies
\[
\mathbf{P}(\varrho[H]_s=\varrho[\widetilde H]_s \mbox{ for a.e. }
s\in[0,t])=1.
\]
\end{definition}

For $\beta\in\Gamma^0_{kl}, Y \in M^0_k$, we write
$J^x(Y, \beta(Y))$ as $J^x(Y, \beta)$. For $\alpha\in\Gamma
^0_{lk}, Z \in M^0_l$,
$J^x(\alpha, Z)$ is defined similarly.

For $k,l\in{\mathbb{N}}$, let
%
%
\begin{eqnarray}\label{50}
V^0_{kl}(x)&=&\inf_{\beta\in\Gamma^0_{kl}}\sup_{Y\in
M^0_k}J^x(Y,\beta),
\\
\label{51}
U^0_{kl}(x)&=&\sup_{\alpha\in\Gamma^0_{lk}}\inf_{Z\in
M^0_l}J^x(\alpha,Z).
\end{eqnarray}
The following shows that the two formulations are equivalent.
\begin{lemma}
\label{lem7}
For every $k,l$, $V^0_{kl}=V_{kl}$ and $U^0_{kl}=U_{kl}$.
\end{lemma}
\begin{pf}
We only show the claim regarding $V_{kl}$.
Let $\beta\in\Gamma_l$. Define $\beta^0\in\Gamma^0_{kl}$ by
letting, for
every $Y\in M^0_k$, $\beta^0[Y]$ be the process component of the
pair $\beta[(Y,k)]$. Clearly, for every $Y\in M^0_{k}$,
$J^x((Y,k),\beta)=J^x(Y,\beta^0)$, whence $\sup_{Y\in
M_k}J^x(Y,\beta)\ge\sup_{Y\in M^0_k}J^x(Y,\beta^0)$, and
$V_{kl}(x)\ge V^0_{kl}(x)$.

Next, let $\beta^0\in\Gamma^0_{kl}$. Define $\beta\dvtx M\to M_l$ as
follows. Given $Y\equiv(Y^0,K)\equiv(A,C,K)\in M$, let
$Y^k=(A,C\wedge
k)$, and set $\beta[Y]=(\beta^0[Y^k],l)$. Note that if, for some
$K$, $Y^0$ and $\widetilde Y^0$ are elements of $M^0_K$ and
$Y^0(s)=\widetilde
Y^0(s)$ on $[0,t]$ then $Y^k(s)=\widetilde Y^k(s)$ on $[0,t]$ and so
$\beta^0[Y^k]_s=\beta^0[\widetilde Y^k]_s$ on $[0,t]$. By definition of
$\beta$, it follows that $\beta\in\Gamma_l$. Also, if $(Y^0,K)\in M_k$
then $K\le k$ and thus $J^x((Y^0,K),\beta)=J^x(Y^0,\beta^0)$. This
shows that $\sup_{Y\in M_k}J^x(Y,\beta)\le\sup_{Y^0\in
M^0_k}J^x(Y^0,\beta^0)$. Consequently, $V_{kl}(x)\le V^0_{kl}(x)$.
\end{pf}

Denote
$V_n=V_{k_nl_n}$ and $U_n=U_{k_nl_n}$.
The following result is proved in Section \ref{sec4}.
\begin{theorem}\label{th3}
For some $n_0 \in\mathbb{N}$,
the family $\{V_n; n \ge n_0\}$ is equicontinuous, and so is the family
$\{U_n; n\ge n_0\}$.
\end{theorem}

Consider the Bellman--Isaacs equations for the upper and,
respectively, lower values of the game with bounded controls, namely
%
%
\begin{eqnarray}\label{08}
&&\cases{ \Lambda_n^+(Du,D^2u)-h=0, &\quad in $G$,\cr
u=g, &\quad on $\partial G$,}
\\
\label{19}
&&\cases{\Lambda_n^-(Du,D^2u)-h=0, &\quad in $G$,\cr
u=g, &\quad on $\partial G$.}
\end{eqnarray}
Solutions to these equations are defined analogously to Definition
\ref{def3}, with $\Lambda^\pm_n$ replacing $\Lambda$, and where
there is no
restriction on the derivatives of the test function, that is,
$\mathcaligr{D}$ is replaced with
${\mathbb{R}}^m\times\mathscr{S}(m)$.
\begin{lemma}\label{lem3}
There exists $n_1 \in{\mathbb{N}}$ such that for each $n \ge n_1$,
$U_n$ is
the unique solution to (\ref{08}), and $V_n$ is the unique solution
to (\ref{19}).
\end{lemma}
\begin{pf}
This follows from a more general result, Theorem \ref{th2} in
Section \ref{sec5}.
\end{pf}
\begin{lemma}\label{lem4}
Any subsequential uniform limit of $U_n$ or $V_n$ is a solution of
(\ref{07}).
\end{lemma}
\begin{pf}
Denote by $U_0$ (resp., $V_0$) a subsequential limit of $U_n$
[$V_n$]. By relabeling, we assume without loss that
$U_n$ (resp., $V_n$) converges to $U_0$ [$V_0$]. We will show that
$U_0$ and $V_0$ are
subsolutions of (\ref{07}). The
proof that these are supersolutions is parallel.

We start with the proof that $U_0$ is a subsolution. Fix $x_0\in G$.
Let $\varphi\in\mathcaligr{C}^2(G)$ be such that $U_0-\varphi$ is strictly
maximized at $x_0$. Assume first that $D\varphi(x_0)\ne0$. Since
$U_n\to U_0$ uniformly, we can find $\{x_n\}\subset G$, $x_n\to
x_0$, where $x_n$ is a local maximum of $U_n-\varphi$ for $n\ge N$. We take
$N$ to be larger than $n_1$ of Lemma \ref{lem3}.
Since by Lemma \ref{lem3} $U_n$ is a subsolution of (\ref{08}), we
have that for $n \ge N$
\[
\Lambda_n^+(D\varphi(x_n),D^2\varphi(x_n))-h(x_n)\le0.
\]
Thus, by Lemma \ref{lem1},
\[
\Lambda(D\varphi(x_0),D^2\varphi(x_0))-h(x_0)\le0
\]
as required.

Next, assume that $D\varphi(x_0)=0$ and $D^2\varphi(x_0)=\lambda I_m$
for some
$\lambda\in{\mathbb{R}}$. In particular,
$\varphi(x)=\varphi(x_0)+\frac\lambda2|x-x_0|^2+o(|x-x_0|^2)$. We
need to show
that
%
%
\begin{equation}
\label{1.1}
-2\lambda-h(x_0)\le0.
\end{equation}
Consider the case $\lambda\ge0$. Fix $\delta>0$ and let
$\psi_\delta(x)=\frac{\lambda+\delta}{2}|x-x_0|^2$. Then $U_0-\psi
_\delta$ has
a strict maximum at $x_0$. Since $U_n\to U_0$ uniformly, we can find
$\{x_n\}\subset G$, $x_n\to x_0$, where $x_n$ is a local maximum of
$U_n-\psi_\delta$. To prove (\ref{1.1}), it suffices to show that for
each $\varepsilon>0$,
%
%
\begin{equation}
\label{1.2}
-2(\lambda+\delta)-\sup_{x\in\mathbb{B}_\varepsilon(x_0)}h(x)\le0.
\end{equation}
To prove (\ref{1.2}), argue by contradiction and assume that it
fails. Then there exists $\varepsilon>0$ such that
%
%
\begin{equation}
\label{1.25}
-2(\lambda+\delta)-\sup_{x\in\mathbb{B}_\varepsilon(x_0)}h(x)>0.
\end{equation}
Let $N \ge n_1$ be such that $|x_n-x_0|<\varepsilon$ for all $n\ge N$.
Since $U_n$
is a subsolution of (\ref{08}),
%
%
\begin{equation}
\label{1.3}
\mu_n:=\Lambda^+_n(D\psi_\delta(x_n),D^2\psi_\delta(x_n))\le h(x_n).
\end{equation}
Also,
%
%
\begin{eqnarray}\label{1.325}\qquad
\mu_n &=& \max_{|b|=1,0\le d\le l_n}\min_{|a|=1,0\le c\le
k_n} \biggl[-\frac12(\lambda+\delta)|a-b|^2\nonumber\\
&&\hspace*{107pt}{} -(\lambda+\delta)(c+d)(a+b)\cdot(x_n-x_0) \biggr]
\\
&\ge&
\min_{|a|=1,0\le c\le k_n} \biggl[-\frac12(\lambda+\delta)
|a-b_n|^2 \biggr] =
-2(\lambda+\delta),\nonumber
\end{eqnarray}
where $b_n=-(x_n-x_0)/|x_n-x_0|$ if $x_n\ne x_0$ and arbitrary
otherwise. Thus by (\ref{1.25}),
%
%
\begin{equation}\label{1.34}
\mu_n > h(x_n).
\end{equation}
However, this contradicts (\ref{1.3}). Hence, (\ref{1.2}) holds and
so (\ref{1.1}) follows.

Consider now the case $\lambda<0$. Let $\delta>0$ be such that
$\lambda+\delta<0$. Let $\psi_\delta$ be as above. Then $U_0-\psi
_\delta$ has
a strict maximum at $x_0$. Fix $\varepsilon>0$. Then one can find
$\gamma<\varepsilon$ such that
%
%
\begin{equation}
\label{1.35}\qquad
U_0(x_0)=U_0(x_0)-\psi_\delta(x_0)>U_0(x)-\psi_\delta(x)\qquad
\forall 0<|x-x_0|\le\gamma.
\end{equation}
Thus, one can find $\eta\in{\mathbb{R}}^m$ such that $0<|\eta
|<\gamma$ and
%
%
\begin{equation}
\label{1.4}
U_0(x_0)>U_0(x)-\psi_\delta(x+\eta)\qquad \forall x\in\partial
\mathbb{B}_\gamma(x_0).
\end{equation}
Let $\psi_{\delta,\eta}(x)=\psi_\delta(x+\eta)$. Let
$x_\eta\in\overline{\mathbb{B}_\gamma(x_0)}$ be a maximum point for
$U_0-\psi_{\delta,\eta}$ over $\overline{\mathbb{B}_\gamma
(x_0)}$. We claim that
%
%
\begin{equation}
\label{1.8}
x_\eta\notin\partial\mathbb{B}_\gamma(x_0) \quad\mbox{and}\quad
x_\eta\ne
x_0-\eta.
\end{equation}
Suppose the claim holds. Then $D\psi_{\delta,\eta}(x_\eta)\ne0$, and
so from the first part of the proof
\[
-2(\lambda+\delta)-h(x_\eta)=\Lambda(D\psi_{\delta,\eta}(x_\eta
),D^2\psi_{\delta,\eta}(x_\eta))-h(x_\eta)\le0.
\]
Since $|x_\eta-x_0|\le\gamma<\varepsilon$, sending $\varepsilon\to
0$ and then
$\delta\to0$ yields (\ref{1.1}).

We now prove (\ref{1.8}). From (\ref{1.4}) and the fact that
$\lambda+\delta<0$,
\[
\sup_{x \in\partial\mathbb{B}_\gamma(x_0)}[U_0(x)-\psi_{\delta
,\eta
}(x)]<U_0(x_0)\le
U_0(x_0)-\psi_{\delta,\eta}(x_0).
\]
Hence, $x_\eta\notin\partial\mathbb{B}_\gamma(x_0)$. Also
\begin{eqnarray*}
U_0(x_0-\eta)-\psi_{\delta,\eta}(x_0-\eta)&=&U_0(x_0-\eta
)<U_0(x_0)+\psi_\delta(x_0-\eta)\\
&\le&
U_0(x_0)-\psi_{\delta,\eta}(x_0),
\end{eqnarray*}
where we used (\ref{1.35}) and the negativity of the functions
$\psi_\delta$ and $\psi_{\delta,\eta}$. This shows that $x_\eta
\ne
x_0-\eta$, and (\ref{1.8}) follows. This completes the proof that
$U_0$ is a subsolution of (\ref{07}).

Finally, the argument for $V_0$ differs only at one point. If we had
$(V_n,\Lambda_n^-)$ instead of $(U_n,\Lambda_n^+)$, then instead of
(\ref{1.325}), we could write
\begin{eqnarray*}
\mu_n
&=&\min_{|a|=1,0\le c\le k_n}\max_{|b|=1,0\le d\le l_n} \biggl[-\frac
12(\lambda+\delta)|a-b|^2\\
&&\hspace*{107pt}{}-(\lambda+\delta)(c+d)(a+b)
\cdot(x_n-x_0) \biggr]
\\
&=&
\max_{|b|=1,0\le d\le
l_n} \biggl[-\frac12(\lambda+\delta)|a_n-b|^2-(\lambda+\delta
)(c_n+d)(a_n+b)\cdot(x_n-x_0) \biggr],
\end{eqnarray*}
where $(a_n,c_n)$ achieves the minimum, and then by choosing
$(b,d)=(-a_n,0)$,
\[
\mu_n\ge-2(\lambda+\delta).
\]
Hence, (\ref{1.34}) is still true. Rest of the argument for the
subsolution property of $V_0$
follows as that for $U_0$.
\end{pf}
\begin{lemma}\label{lem5}
Fix $x\in\overline{G}$.

\begin{longlist}
\item One can choose $(k_n,l_n)$ in such a way that
$\limsup_{n\to\infty}V_n(x)\le V(x)$.
\item One can choose $(k_n,l_n)$ in such a way that
$\liminf_{n\to\infty}V_n(x)\ge V(x)$.

Similar statements hold for $U_n(x)$ and $U(x)$.
\end{longlist}
\end{lemma}
\begin{pf}
We prove (i) and (ii). The statements regarding $U_n(x)$ and
$U(x)$ are proved analogously.

\begin{longlist}
\item
Fix $k$. Since $\Gamma= \bigcup_{l\ge1}\Gamma_l$, we have that
given $\varepsilon$,
\begin{eqnarray*}
V(x)&\ge&\inf_{\beta\in\Gamma}\sup_{Y\in M_k}J^x(Y,\beta)\\
&\ge&\inf_{\beta\in\Gamma_l}\sup_{Y\in M_k}J^x(Y,\beta
)-\varepsilon\\
&=&V_{kl}(x)-\varepsilon,
\end{eqnarray*}
for all $l$ sufficiently large. This shows
$V(x)\ge\limsup_{l\to\infty}V_{kl}(x)$, and (i) follows.
\item Fix $\varepsilon$. For each $(k,l)\in{\mathbb{N}}^2$, let
$\beta_{kl}\in\Gamma
_l$ be
such that
%
%
\begin{equation}\label{22}
\sup_{Y\in M_k}J^x(Y,\beta_{kl})\le\inf_{\beta\in\Gamma_{l}}\sup
_{Y\in
M_k}J^x(Y,\beta)+\varepsilon.
\end{equation}
Fix $l$. Let $\beta_l$ be defined by
\[
\beta_l[Y]=\beta_{kl}[Y],\qquad  Y\in M_k\setminus M_{k-1},  k\in
{\mathbb{N}},
\]
where we define $M_0$ to be the empty set. Then $\beta_l\in\Gamma_{l}$.
Since $M = \bigcup_{k \ge1}M_k$, we have that
the following holds provided that $k$
is sufficiently large
\begin{eqnarray*}
V(x)&\le&\sup_{Y\in M}J^x(Y,\beta_l)\\
&\le&\sup_{Y\in M_k}J^x(Y,\beta_l)+\varepsilon\\
&=&\max_{j\le k}\sup_{Y\in M_j\setminus M_{j-1}}J^x(Y,\beta
_{jl})+\varepsilon
\\
&\le&\inf_{\beta\in\Gamma_{l}}\sup_{Y\in M_k}J^x(Y,\beta
)+2\varepsilon,
\end{eqnarray*}
where the last inequality follows from (\ref{22}). This shows that,
for every $l$, $V(x)\le\liminf_kV_{kl}(x)$. The result follows.\qed
\end{longlist}
\noqed\end{pf}
\begin{pf*}{Proof of Theorem \protect\ref{th1}}
The statement that $U$ and $V$ are solutions of (\ref{07}) follows from
Theorem \ref{th3}, Lemmas \ref{lem3}, \ref{lem4} and uniqueness of
solutions of (\ref{07}), established in \cite{pssw}. The latter result
also yields $U=V$.
\end{pf*}

\section{Equicontinuity}\label{sec4}

In this section, we prove Theorem \ref{th3}.
With an eye toward estimates needed in Section
\ref{sec5} we will consider a somewhat more general setting.
Thanks to Lemma \ref{lem7} we may, and will use the value functions
(\ref{50}), (\ref{51}), defined using
simple controls and strategies (Definition \ref{def2}).
Given $X$ defined as in (\ref{01})
for some $Y, Z \in M^0$, we let for $\gamma\in[0, 1)$, $X^{\gamma} =
X + \gamma\widetilde W$. Define $\tau^{\gamma}$ and $J_{\gamma}$ as below
(\ref{23}) but with $X$ replaced with $X^{\gamma}$. Also denote by
$U^{\gamma}_{kl}$ and $V^{\gamma}_{kl}$ the expressions in (\ref{50}),
(\ref{51})
with $J$ replaced with $J_{\gamma}$. We write $U_n^{\gamma}
= U^{\gamma}_{k_nl_n}$, $V_n^{\gamma} = V^{\gamma}_{k_nl_n}$. Theorem
\ref{th3} is an immediate consequence of the following more general
result.
\begin{theorem}\label{th3plus}
For some $n_2 \in\mathbb{N}$,
the family $\{V_n^{\gamma}, U_n^{\gamma}; n \ge n_2, \gamma\in[0,
1)\}$ is equicontinuous.
\end{theorem}

In what follows, we will suppress $\gamma$ from the notation unless
there is a scope for confusion. We start by showing
that the value functions are uniformly\vspace*{1pt} bounded. To this end,
fix $a^0\in S^{m-1}$, and note that the constant process $Y^0:=(a^0,1)$
is in
$M^0$.
\begin{lemma}
\label{lem01}
There exists a constant $c_1<\infty$ such that
\[
\mathbf{E}[\tau(x,Y^0,Z)^2]\le c_1, \qquad x\in\overline G,  Z\in M^0,
\gamma
\in[0,1).
\]
\end{lemma}
\begin{pf}
We only present the proof for the case $\gamma= 0$. The
general case follows upon minor modifications. Denote by $m_0$
the diameter of $G$. Fix $T > m_0$. By~(\ref{01}), with
$\alpha_t=a^0\cdot B_t$, on the event $\tau>T$ one has
\[
\int_0^T(1-\alpha_s)\,dW_s+\int_0^T(1+\alpha_s)\,ds\le a^0\cdot
(X_T-X_0)<m_0.
\]
Consider the $\{\mathcaligr{F}_t\}$-martingale, $M_t = \int_0^t
(1-\alpha_s)\,
dW_s$, with $\langle M\rangle_t = \int_0^t (1-\alpha_s)^2 \,ds$.
On the event $\langle M \rangle_T < T$,
\[
\int_0^T(1 + \alpha_s) \,ds = 2T - \int_0^T (1-\alpha_s)\,ds \ge T.
\]
So on the set $\{\langle M \rangle_T < T; \tau> T\}$ we have $|M_T|
\ge
T - m_0$. Letting $\sigma= \inf\{s\dvtx\langle M \rangle_s \ge T\}$,
%
%
\begin{eqnarray}\label{ab442}
\mathbf{P}(\tau> T; \langle M \rangle_T < T) &\le&
\mathbf{P}(|M_{T\wedge\sigma}| \ge T - m_0) \nonumber\\[-8pt]\\[-8pt]
&\le&\frac{m_1\mathbf{E}
\langle
M \rangle^2_{T\wedge\sigma}}{(T-m_0)^4} \le\frac{m_1
T^2}{(T-m_0)^4}.\nonumber
\end{eqnarray}
We now consider the event $\{\tau> T ; \langle M \rangle_T \ge T\}$. One
can find $m_2, m_3 \in(0, \infty)$ such that for all
nondecreasing, nonnegative processes $\{\widehat\gamma_t\}$,
%
%
\begin{equation}\label{ab455}
\mathbf{P}\bigl(H_s + \widehat\gamma_s \in(-m_0, m_0); 0
\le s \le T \bigr) \le m_2e^{-m_3 T},
\end{equation}
where $H$ is a one-dimensional Brownian motion. Letting $\gamma_t =
\int_0^t (1+D_s) (1 + \alpha_s) \,ds$, where $Z = (B,D)$, we see that
\[
\{\tau> T; \langle M \rangle_T \ge T\} \subset
\{M_s + \gamma_s \in(-m_0, m_0), 0 \le s \le T; \langle M \rangle
_T \ge T\}.
\]
For $u\ge0$, let $S_u = \inf\{s\dvtx\langle M \rangle_s > u\}$. Then, with
$\widehat\gamma_s = \gamma_{S_s}$,
\[
\mathbf{P}(\tau> T; \langle M \rangle_T \ge T) \le
\mathbf{P}\bigl(H_s + \widehat\gamma_s \in(-m_0, m_0); 0 \le s \le T\bigr) \le
m_2e^{-m_3 T},
\]
where the last inequality follows from (\ref{ab455}). The result now
follows on combining the above display with (\ref{ab442})
\end{pf}

The inequality $J(x, Y^0, Z) \le|h|_{\infty} \mathbf{E}(\tau(x,
Y^0, Z)) +
|g|_{\infty}$,
where $|h|_{\infty} = {\sup_x} |h(x)|$ and $|g|_{\infty} = {\sup_x} |g(x)|$,
immediately implies the following.
\begin{corollary}
\label{cor01}
There exists a constant $c_2<\infty$ such that $|V_n^{\gamma}(x)|\vee
|U_n^{\gamma}(x)|\le c_2$, for
all $x\in\overline G$, $\gamma\in[0, 1)$ and $n\in{\mathbb{N}}$.
\end{corollary}

The idea of the proof of equicontinuity, explained in a heuristic
manner, is as follows. Let $x_1$ and $x_2$ be in $G$, let
$\varepsilon=|x_1-x_2|$, and let $\delta>0$. Consider the game with bounded
controls for which $V_n$ is the lower value function, for some
$n\in{\mathbb{N}}$. Let the minimizing player select a strategy
$\beta^n$ that
is $\delta$-optimal for the initial position $x_1$; namely $\sup_{Y
\in M^0_{k_n}} J^{x_1}(Y, \beta^n) \le V_n(x_1) + \delta$. Denote
the exit time by $\tau_1=\tau^{x_1}(Y,\beta^n)$ and the exit
position by $\xi_1=X^{x_1}(\tau_1)$. Now, modify the strategy is
such a way that the resulting control $Z=\beta^n[Y]$ is only
affected for times $t\ge\tau_1$. This way, the payoff incurred
remains unchanged. Thus, denoting the modified strategy by
$\widetilde\beta^n$, we have, for every $Y\in M^0_{k_n}$,
\[
J^{x_1}(Y,\widetilde\beta^n)\le V_n(x_1)+\delta.
\]
Given a
point $\xi_2$ located inside $G$, $\varepsilon$ away from $\xi_1$,
and a
new state process which, at time $\tau_1$ is located at $\xi_2$, the
modified strategy attempts to force this process to exit the domain
soon after $\tau_1$ and with a small displacement from $\xi_2$
(provided that $\varepsilon$ is small).

Let now the maximizing player select a control $Y^n$ that is
$\delta$-optimal for playing against $\widetilde\beta^n$, when
starting from
$x_2$. This control is modified after the exit time $\tau^{x_2}(Y^n,
\widetilde\beta^n)$
in a similar
manner to the above. Denoting the modified control by~$\widetilde Y^n$, we
have
\[
V_n(x_2)\le J^{x_2}(\widetilde Y^n,\widetilde\beta^n)+\delta.
\]
Hence, $V_n(x_2)-V_n(x_1)\le J^{x_2}(\widetilde
Y^n,\widetilde\beta^n)-J^{x_1}(\widetilde Y^n,\widetilde\beta
^n)+2\delta$. One can thus
estimate the modulus of continuity of $V_n$ by analyzing the payoff
incurred when $(\widetilde Y^n,\widetilde\beta^n)$ is played, considering
simultaneously two state processes, starting from $x_1$ and~$x_2$.
The form (\ref{01}) of the dynamics ensures that the processes
remain at relative position $x_1-x_2$ until, at time $\sigma$, one of
them leaves the domain. The difference between the running payoffs
incurred up to that time can be estimated in terms of~$\varepsilon$, the
modulus of continuity of $h$, and the expectation of $\sigma$. It is not
hard to see that the latter is uniformly bounded, owing to Corollary
\ref{cor01} and the boundedness of $h$ away from zero. By
construction, one of the players will now attempt to force the state
process that is still in $G$ to exit. If one can ensure that exit
occurs soon after $\sigma$ and with a small displacement (uniformly in
$n$), then the running payoff incurred between time $\sigma$ and the
exit time is small, and the difference between the terminal payoffs
is bounded in terms of $\varepsilon$ and the modulus of continuity of $g$,
resulting in an estimate that is uniform in $n$.

This argument is made precise in the proof of the theorem. Lemmas
\ref{lem02} and \ref{lem03} provide the main tools for showing that
starting at a state near the boundary, each player may force exit
within a short time and with a small displacement. To state these
lemmas, we first need to introduce some notation.

We have assumed that $G$ is a bounded $C^2$ domain in ${\mathbb
{R}}^m$. Thus,
there exist $\overline\rho\in(0,\frac18)$, $k\in{\mathbb{N}}$, $z_j\in
\partial G$, $E_j\in\mathcaligr{O}(m)$, $\xi_j\in C^2({\mathbb{R}}^{m-1})$,
$j=1,\ldots,k$, such
that, with $\mathbb{B}_j=\mathbb{B}_{\overline\rho}(z_j)$, $j=1,\ldots
,k$, one has
$\partial G\subset\bigcup_{j=1}^k\mathbb{B}_j$, and
\[
G\cap\mathbb{B}_j=\{E_jy\dvtx y_1>\xi_j(y_2,\ldots,y_m)\}\cap
\mathbb{B}_j,\qquad  j=1,\ldots,k.
\]
Here, $\mathcaligr{O}(m)$ is the space of $m\times m$ orthonormal matrices.
Define for $j=1,\ldots,k$, $\widetilde\varphi_j\dvtx{\mathbb{R}}^m\to
{\mathbb{R}}$ as
\[
\widetilde\varphi_j(y)=y_1-\xi_j(y_2,\ldots,y_m),\qquad  y\in
{\mathbb{R}}^m.
\]
Let $\varphi_j(x)=\widetilde\varphi_j(E_j^{-1}x)$, $x\in{\mathbb
{R}}^m$. Then
$|D\varphi_j(x)|\ge1$, $x\in{\mathbb{R}}^m$, $j=1,\ldots,k$. Furthermore,
\[
G\cap\mathbb{B}_j=\{x\dvtx\varphi_j(x)>0\}\cap\mathbb{B}_j,\qquad
j=1,\ldots,k.
\]
Let $0<\rho_0<\overline\rho$ be such that $\partial
G\subset\bigcup_{j=1}^k\mathbb{B}_{\rho_0}(z_j)$. For $\varepsilon
>0$, denote
\[
\mathbf{X}_\varepsilon=\{(x_1,x_2)\dvtx x_1\in\partial G,  x_2\in G,
|x_1-x_2|\le\varepsilon
\}.
\]
Let $\underline j\dvtx \partial G\to\{1,\ldots,k\}$ be a measurable map
with the property
\[
x\in\mathbb{B}_{\rho_0}\bigl(z_{\underline j(x)}\bigr)\qquad  \mbox{for all }
x\in\partial G.
\]
For existence of such a map see, for example, Theorem 10.1 of
\cite{EK}. Then, for every
$\varepsilon\le\rho_1:=\frac{\overline\rho-\rho_0}{4}$,
%
%
\begin{equation}\label{34}
(x_1,x_2)\in\mathbf{X}_\varepsilon\qquad\mbox{implies }
\overline{\mathbb{B}_{\rho_1}(x_i)}\subset\mathbb{B}_{\underline
j(x_1)},\qquad i=1,2.
\end{equation}
For $j=1,\ldots,k$ and $x_0\in\mathbb{B}_j$, define
$\psi^{x_0}_j\dvtx{\mathbb{R}}^m\to{\mathbb{R}}$ as
\[
\psi^{x_0}_j(x)=\varphi_j(x)+|x-x_0|^2.
\]
Also, note that $|D\psi^{x_0}_j|\ge\frac12$ in $\mathbb{B}_j$. Define
$\pi_j^{x_0}\dvtx{\mathbb{R}}^m\to\mathcaligr{S}^{m-1}$ such that it is
Lipschitz, and
%
%
\begin{equation}\label{31}
\pi_j^{x_0}(x)=-\frac{D\psi^{x_0}_j(x)}{|D\psi^{x_0}_j(x)|},\qquad
x\in\mathbb{B}_j.
\end{equation}

Given a strategy $\beta$, and a point $x_2$, we seek a control
$Y=(A,C)$ that forces a state process starting from $x_2$ to exit in
a short time and with a small displacement from $x_2$ (provided that
$x_2$ is close to the boundary). We would like to determine $Y$ via
the functions $\pi_j$ just constructed, in such a way that the
following relation holds:
%
%
\begin{equation}\label{39}
A(t)=\pi(X(t)), \qquad C(t)=c^0,
\end{equation}
where $\pi=\pi^{x_2}_{\underline j(x_1)}$ and $c^0>0$ is some constant.
Making $A$ be oriented in the negative direction of the gradient of
$\varphi_j$ allows us to show that the state is ``pushed'' toward the
boundary. The inclusion of a quadratic term in $\psi_j^{x_0}$
ensures in addition that the sublevel sets $\{\psi_j<a\}$ are
contained in a small vicinity of $x_0$, provided smallness of $a$
and $\operatorname{dist}(x_0,\partial G)$. The latter property enables us
to show
that the process does not wander a long way along the boundary
before exiting.

The difficulty we encounter is that due to the feedback nature of
$A$ in (\ref{39}) we cannot ensure (local) existence of solutions of
the set of (\ref{01}) and (\ref{39}). Take, for example, a
strategy $\beta$ that is given as $\beta[Y]_t = b(Y_t)$, $Y \in
M^0$, where $b$ is some measurable map from $\mathcaligr{H}$ to $\mathcaligr{H}$.
Along with (\ref{01}), and (\ref{39}) this
defines $X$ as a solution to an SDE with general measurable
coefficients. However,
as is well known, the SDE may not admit any solution in this
generality. To overcome this problem, we will construct a $Y$ that approximates
the $Y$ we seek in (\ref{39}) via a time
discretization.

Let $\bolds{\Sigma}_\varepsilon$ denote the collection of all quintuples
$\bolds{\sigma}=(\sigma,\xi_1,\xi_2,\beta,Y)$ such that $\sigma$ is
an a.s. finite
$\mathcaligr{F}_t$ stopping time, $\xi_1$ and $\xi_2$ are
$\mathcaligr{F}_\sigma$-measurable random variables satisfying $(\xi
_1,\xi_2)\in
\mathbf{X}_\varepsilon$ a.s., $\beta\in\Gamma^0$, and $Y\in M^0$.
Fix $\gamma
\in
[0, 1)$.

Let $\bolds{\sigma}=(\sigma,\xi_1,\xi_2,\beta,Y)\in\bolds
{\Sigma}
_{\rho_1}$ be
given. Let $j^*(\omega)=\underline j(\xi_1(\omega))$. Denote
\[
\Phi\equiv\Phi(\omega,\cdot)=\varphi
_{j^*(\omega)},\qquad
\Psi\equiv\Psi(\omega,\cdot)=\psi^{\xi_2(\omega)}_{j^*(\omega
)},\qquad
\Pi\equiv\Pi(\omega,\cdot)=\pi^{\xi_2(\omega)}_{j^*(\omega)}.
\]
We define a sequence of processes $(X^{(i)},Y^{(i)})_{i\ge0}$ as
follows. Let $Y_t^{(0)}\equiv(A_t^{(0)},\break C_t^{(0)})$ be given by
\[
Y_t^{(0)} = \cases{Y_t, &\quad $t<\sigma$,\cr
(\Pi(\xi_2),c^0), &\quad $t\ge\sigma$.}
\]
The constant $c^0$ above will be chosen later. Denote
$(B^{(0)},D^{(0)})=\beta[Y^{(0)}]$. Define, for $t\ge\sigma$,
\begin{eqnarray*}
X_t^{(0)} &=& \xi_2+\int\mathbf{1}_{[\sigma,t]}(s)
\bigl(\bigl[A^{(0)}_s-B^{(0)}_s\bigr]\,dW_s
+ \gamma d\widetilde W_s \bigr)
\\
&&{} +\int_\sigma^t\bigl[C^{(0)}_s+D^{(0)}_s\bigr]\bigl[A^{(0)}_s+B^{(0)}_s\bigr]\,ds.
\end{eqnarray*}
The process $X^{(0)}$ can be defined arbitrarily for $t<\sigma$. Set
$\eta_0=\sigma$. We now define recursively, for all $i\ge1$,
%
\begin{eqnarray}\hspace*{33pt}
\label{33}
\eta_i &=&
(\eta_{i-1}+\varepsilon)
\nonumber\\[-8pt]\\[-8pt]
&&{}\wedge\inf \biggl\{
t\ge\eta_{i-1}\dvtx
\bigl|X^{(i-1)}_t-X^{(i-1)}_{\eta_{i-1}}\bigr|\ge\varepsilon \mbox
{ or }
\int_{\eta_{i-1}}^t R_s^{(i)}\,ds \ge t-\eta_{i-1} \biggr\},\nonumber
\end{eqnarray}
where
$R^{(i)}_s=[C^{(i-1)}_s+D^{(i-1)}_s]D\Psi(X^{(i-1)}_s)\cdot
[A^{(i-1)}_s-\Pi(X_s^{(i-1)})]$,
\begin{eqnarray}\label{32}
X^{(i)}_t &=&
X^{(i-1)}_t, \qquad Y^{(i)}_t=Y^{(i-1)}_t, \qquad 0\le
t<\eta_i,\nonumber\\
Y^{(i)}_t
&\equiv&
\bigl(A^{(i)}_t,C^{(i)}_t\bigr)=\bigl(\Pi\bigl(X^{(i-1)}_{\eta_i}\bigr),c^0\bigr),\qquad
t\ge\eta_i,\nonumber\\
\bigl(B^{(i)},D^{(i)}\bigr) &=& \beta\bigl[Y^{(i)}\bigr],
\\
X^{(i)}_t &=&
X^{(i-1)}_{\eta_i}+\int\mathbf{1}_{[\eta_i,t]}(s)
\bigl(\bigl[A_s^{(i)}-B_s^{(i)}\bigr]\,dW_s+ \gamma d\widetilde W_s \bigr)\nonumber\\
&&{} +\int_{\eta_i}^t\bigl[A^{(i)}_s+B^{(i)}_s\bigr]\bigl[C^{(i)}_s+D^{(i)}_s\bigr]\,ds,\qquad
t>\eta_i.\nonumber
\end{eqnarray}
It is easy to check that $\eta_0<\eta_1<\cdots$ and $\eta_i\to
\infty$
a.s. Define $X_t=X^{(i)}_t$, $Y_t=Y^{(i)}_t$ if $t\le\eta_i$. Let
$\rho=\rho_1^2$ and
%
%
\begin{eqnarray}\label{4.1}
\tau_{\rho} &=& \inf\{t\ge\sigma\dvtx\Psi(X_t)\ge\rho\},\nonumber\\
\tau_G &=& \inf\{t\ge\sigma\dvtx X_t\in G^c\},\\
\tau &=& \tau_{\rho}\wedge\tau_G.\nonumber
\end{eqnarray}
Define
\begin{eqnarray*}
\overline Y_t &\equiv& (\overline A_t,\overline C_t)=
\cases{Y_t, &\quad $t<\tau$, \cr
(a^0,c^0), &\quad $t\ge\tau$,}
\\
(\overline B,\overline D)&=&\beta(\overline Y),
\end{eqnarray*}
where $a^0$ is as fixed at the beginning of the section. Let
\[
\overline X_t=\xi_2+\int\mathbf{1}_{[\sigma,t]}(s) ([\overline A_s-\overline
B_s]\,dW_s+ \gamma \,d\widetilde W_s ) +\int_\sigma^t[\overline C_s+\overline
D_s][\overline
A_s+\overline B_s]\,ds.
\]
We write
\[
\overline X=\overline X[\sigma,\xi_1,\xi_2,\beta,Y],\qquad
\overline Y=\overline Y[\sigma,\xi_1,\xi_2,\beta,Y].
\]
Note that if $\overline\tau_{\rho}$ and $\overline\tau_G$ are defined by
(\ref{4.1}) upon replacing $X$ by $\overline X$ then
$\overline\tau:=\overline\tau_{\rho}\wedge\overline\tau_G=\tau_{\rho}\wedge
\tau
_G=\tau$,
because $\overline X$ differs from $X$ only after time $\tau$. We write
$\overline\tau=\overline\tau[\sigma,\xi_1,\xi_2,\beta,Y]$. Similar
notation will
be used for $\overline\tau_{\rho}$ and $\overline\tau_G$.
\begin{lemma}
\label{lem02}
There exists a $c^0 \in(0, \infty)$ and
a modulus $\vartheta$ such that for every $\varepsilon\in(0, \rho
_1)$, and
$\gamma\in[0, 1)$, if $\bolds{\sigma}=(\sigma,\xi_1,\xi_2,\beta
,Y)\in\bolds{\Sigma}_\varepsilon$
and $\overline\tau_G=\overline\tau_G[\bolds{\sigma}]$, one has:

\begin{longlist}
\item
$\mathbf{E}\{\overline\tau_G-\sigma |\mathcaligr{F}_\sigma\}\le
\vartheta(\varepsilon)$,
\item
$\mathbf{E}\{|\overline X-\xi_2|^2_{*,\overline\tau_G} |\mathcaligr
{F}_\sigma\}\le\vartheta(\varepsilon)$,
where $|\overline X-\xi_2|_{*,\overline\tau_G} = {\sup_{t \in[\sigma, \overline
\tau_G]}}|\overline
X (t) - \xi_2|$.
\end{longlist}
\end{lemma}

Proof of the lemma is provided after the proof of Theorem \ref{th3plus}.

Next, we construct a strategy $\beta^*\in\Gamma^0$ with analogous
properties. Here, existence of solutions is not an issue, and
discretization is not needed.

Fix $(x_1,x_2)\in\mathbf{X}_{\rho_1}$. Let $j=\underline j(x_1)$,
$\widetilde\Psi=\psi_j^{x_2}$, and $\widetilde\Pi=\pi_j^{x_2}$. Given
$Y=(A,C)\in
M^0$, let $\widetilde X$ solve
\[
\widetilde X_t=x_2+\int_0^t \bigl([A_s-\widetilde\Pi(\widetilde
X_s)]\,dW_s + \gamma \,d\widetilde W_s \bigr)+\int
_0^t[C_s+d^0][A_s+\widetilde\Pi
(\widetilde X_s)]\,ds,
\]
where $d^0$ is a constant to be determined later. Let
%
%
\begin{eqnarray}\label{5.1}
\widetilde\tau_{\rho}
&=&\inf\{t\dvtx\widetilde\Psi(\widetilde X_t)\ge\rho\},\nonumber\\
\widetilde\tau_G
&=&\inf\{t\dvtx\widetilde X_t\in G^c\},
\\
\widetilde\tau
&=&\widetilde\tau_{\rho}\wedge\widetilde\tau_G.\nonumber
\end{eqnarray}
Define $Z^*\in M^0$ as
\[
Z^*_s\equiv(B^*_s,D^*_s)=
\cases{(\widetilde\Pi(\widetilde X_s),d^0), &\quad $s<\widetilde\tau
$,\vspace*{2pt}\cr
(a^0, d^0), &\quad $s\ge\widetilde\tau$.}
\]
Note that $\beta^*[Y](s):=\widetilde Z_s$, $s\ge0$ defines a
strategy. Let
\[
X^*_t=x_2+\int_0^t ([A_s-B^*_s]\,dW_s + \gamma \,d\widetilde W_s
)+\int
_0^t[C_s+d^0][A_s+B^*_s]\,ds.
\]
Define $\tau_{\rho}^*$, $\tau_{G}^*$ and $\tau^*$
by replacing $\widetilde X$ with $X^*$ in (\ref{5.1}), and note that
$\tau^*=\widetilde\tau$. To make the dependence explicit, we write
\[
X^*=X^*[x_1,x_2,Y,\overline W], \qquad  Z^*=Z^*[x_1,x_2,Y,\overline
W],\qquad
\tau^*=\tau^*[x_1,x_2,Y,W].
\]
\begin{lemma}
\label{lem03}
There exists $d^0 \in(0, \infty)$ and a modulus $\widetilde\vartheta
$ such
that for all $\varepsilon\in(0,
\rho_1)$, $Y\in M^0$, $(x_1,x_2)\in\mathbf{X}_\varepsilon$, if
$\tau^*_G=\tau^*_G[x_1,x_2,Y,W]$, one has:

\begin{longlist}
\item
$\mathbf{E}[\tau^*_G]\le\widetilde\vartheta(\varepsilon)$,
\item $\mathbf{E}[|X^*-x_2|^2_{*,\tau^*_G} ]\le\widetilde\vartheta
(\varepsilon)$,
where $|X^*-x_2|_{*,\tau^*_G} = {\sup_{t \in[0, \tau^*_G]}}|X^*(t) - x_2|$.
\end{longlist}
\end{lemma}

The proof of Lemma \ref{lem03} is very similar to (in fact somewhat
simpler than) the proof of Lemma \ref{lem02}, and therefore will be
omitted.

%
%
\begin{figure}[b]

\includegraphics{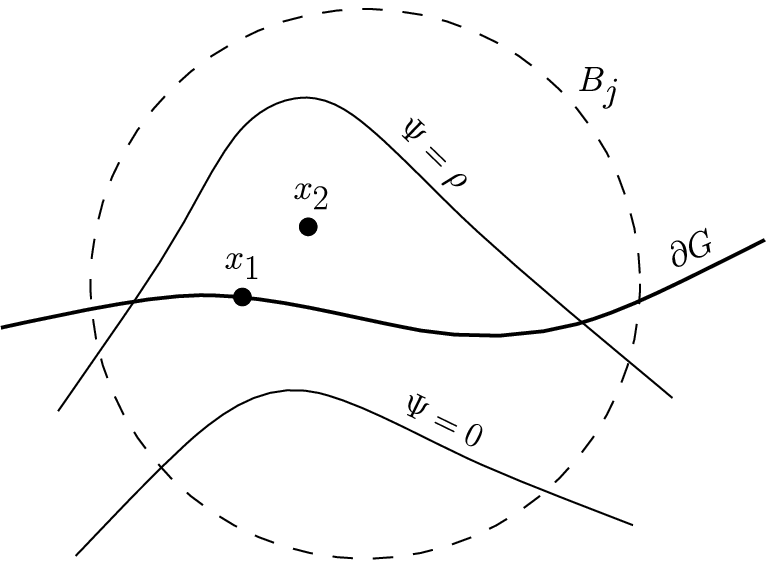}

\caption{}
\label{figure20}
\end{figure}

If $\sigma$ is an a.s. finite $\{\mathcaligr{F}_t\}$-stopping time and
$(\xi_1,\xi_2)$ are $\mathcaligr{F}_\sigma$-measurable random
variables such that
$(\xi_1,\xi_2)\in\mathbf{X}_{\rho_1}$ a.s., then we define the
$\mathcaligr{G}_t=\mathcaligr{F}_{t+\sigma}$ adapted processes
\begin{eqnarray*}
\overline X{}^*_t &=& X^*[\xi_1,\xi_2,\widehat Y_\sigma,\widehat W_\sigma](t),
\\
\overline Z{}^*_t &=& Z^*[\xi_1,\xi_2,\widehat Y_\sigma,\widehat W_\sigma](t),
\end{eqnarray*}
where $\widehat Y_\sigma(t)=Y(t+\sigma)$ and $\widehat W_\sigma
(t)=\overline W(t+\sigma)-\overline W(\sigma)$,
$t\ge0$. To make the dependence explicit, write
\[
\overline X{}^*=\overline X{}^*[\sigma,\xi_1,\xi_2,Y],\qquad
\overline Z{}^*=\overline Z{}^*[\sigma,\xi_1,\xi_2,Y].
\]
\begin{pf*}{Proof of Theorem \protect\ref{th3}}
Fix $x_1,x_2\in G$ and $\gamma\in[0,1)$. We will suppress $\gamma$
from the notation. Assume that $|x_1-x_2|=\varepsilon<\rho_1$, so
that Lemmas
\ref{lem02} and \ref{lem03} are in force (see Figure \ref{figure20}). Let $n_0$ be large enough so
that $l_n, k_n \ge\max\{c^0, d^0\}$ for
all $n \ge n_0$.
Given
$\delta\in(0,1)$ and $n\ge n_0$, let $\beta^n\in\Gamma^0_{k_nl_n}$ be
such that
\[
\sup_{Y\in M^0_{k_n}}J^{x_1}(Y,\beta^n)-\delta\le V_n(x_1)\le c_1.
\]
For $Y\in M^0$ write $\tau^{1,n}(Y):=\tau^{x_1}(Y,\beta^n)$ and
$X^{1,n}(Y):=X^{x_1}(Y,\beta^n)$. Note that
\[
\underline h \mathbf{E}[\tau^{1,n}(Y)]-|g|_\infty\le c_1+1,
\]
hence for every $n$ and every $Y\in M^0_{k_n}$,
%
%
\begin{equation}\label{35}
\mathbf{E}[\tau^{1,n}(Y)]\le m_1,
\end{equation}
where $m_1<\infty$ is a constant that does not depend on $n$.

Define $\widetilde\beta^n\in\Gamma^0_{k_nl_n}$ as follows.
For $Y\in M^0$, let
\begin{eqnarray*}
\xi_1^{1,n}(Y)&=&X^{1,n}_{\tau^{1,n}(Y)}(Y),\qquad
\xi_2^{1,n}(Y)=\xi_1^{1,n}(Y)+x_2-x_1,
\\
\widetilde\beta^n[Y]_t &=&
\cases{\beta^n[Y]_t, &\quad $t<\tau^{1,n}(Y)$,\vspace*{2pt}\cr
\overline Z{}^*[\tau^{1,n}(Y),\xi_1^{1,n}(Y),\xi_2^{1,n}(Y),Y],
&\quad $t\ge\tau^{1,n}(Y),  \xi_2^{1,n}(Y)\in\overline G$,\cr
\mbox{arbitrarily defined}, &\quad $t\ge\tau^{1,n}(Y), \xi
_2^{1,n}(Y)\in\overline
G^c$.}
\end{eqnarray*}
Note that for every $Y\in M^0_{k_n}$,
$J^{x_1}(Y,\beta^n)=J^{x_1}(Y,\widetilde\beta^n)$. Next, choose
$Y^n\in
M^0_{k_n}$ such that
\[
V_n(x_2)\le\sup_{Y\in M^0_{k_n}}J^{x_2}(Y,\widetilde\beta^n)
\le J^{x_2}(Y^n,\widetilde\beta^n)+\delta.
\]
Let $\tau^{2,n}=\tau^{x_2}(Y^n,\widetilde\beta^n)$, and
$X^{2,n}=X^{x_2}(Y^n,\widetilde\beta^n)$. Let
\[
\xi_2^{2,n}=X^{2,n}_{\tau^{2,n}},\qquad
\xi_1^{2,n}=\xi_2^{2,n}+x_1-x_2.
\]
Define $\widetilde Y^n\in M^0_{k_n}$ as
\[
\widetilde Y^n_t=\cases{Y^n_t, &\quad $t<\tau^{2,n}$,\vspace*{2pt}\cr
\overline Y[\tau^{2,n},\xi_2^{2,n},\xi_1^{2,n},\widetilde\beta
^n,Y^n](t), &\quad
$t\ge\tau^{2,n},  \xi_1^{2,n}\in\overline G$,\vspace*{2pt}\cr
\mbox{arbitrarily defined}, &\quad $t\ge\tau^{2,n},  \xi_1^{2,n}\in
\overline
G^c$.}
\]
Note that $J^{x_2}(Y^n,\widetilde\beta^n)=J^{x_2}(\widetilde
Y^n,\widetilde\beta^n)$.
Thus
%
%
\begin{equation}
\label{8.1}
V_n(x_2)-V_n(x_1)-2\delta\le
J^{x_2}(\widetilde Y^n,\widetilde\beta^n)- J^{x_1}(\widetilde
Y^n,\widetilde\beta^n).
\end{equation}
For $k=1,2$, let
\[
\sigma^{k,n}=\tau^{x_k}(\widetilde Y^n,\widetilde\beta^n),\qquad
\widetilde X^{k,n}=X^{x_k}(\widetilde Y^n,\widetilde\beta^n),\qquad
\Xi^{k,n}=\widetilde X^{k,n}_{\sigma^{k,n}}.
\]
For $m_0\ge0$, let $\vartheta_g(m_0)=\sup\{|g(x)-g(y)|\dvtx x,y\in
\partial G,|x-y|\le
m_0\}$ and $\vartheta_h(m_0)=\sup\{|h(x)-h(y)|\dvtx x,y\in G,|x-y|\le
m_0\}$. Using
(\ref{35}), the right-hand side of (\ref{8.1}) can be bounded by
%
%
\begin{equation}\hspace*{33pt}
\label{8.2}
\mathbf{E}\vartheta_g(|\Xi^{1,n}-\Xi^{2,n}|)+c_3\vartheta
_h(\varepsilon)+ |h|_{\infty}
\mathbf{E}[(\sigma^{1,n}\vee\sigma^{2,n})-(\sigma^{1,n}\wedge
\sigma^{2,n})].
\end{equation}
On the set $\sigma^{1,n}\le\sigma^{2,n}$, we have $|\Xi^{1,n}-\Xi^{2,n}|
\le\varepsilon+|\widetilde X_{\sigma^{1,n}}^{2,n}-\widetilde
X_{\sigma^{2,n}}^{2,n}|$.
Hence, by Lem\-ma \ref{lem03}(ii),
\[
\mathbf{E} \bigl[|\Xi^{1,n}-\Xi^{2,n}|^2 \mathbf{1}_{\{\sigma
^{1,n}\le\sigma^{2,n}\}
} \bigr]\le\vartheta_1(\varepsilon)
\]
for some modulus $\vartheta_1$. Using Lemma \ref{lem02}(ii), a similar
estimate holds on the complement set, and consequently, the first term
of (\ref{8.2}) is bounded by
$\vartheta_2(\varepsilon)$, for some modulus $\vartheta_2$. By
Lemmas \ref{lem02}(i)
and \ref{lem03}(i), the last term of (\ref{8.2}) is bounded by
$|h|_{\infty}(\vartheta(\varepsilon)+\widetilde\vartheta
(\varepsilon))$. Hence,
$V_n(x_2)-V_n(x_1)\le2\delta+\vartheta_3(|x_1-x_2|)$ for some modulus
$\vartheta_3$, and the equicontinuity of $\{V_n^{\gamma}; n,\gamma\}$
follows on sending
$\delta\to0$. The proof of equicontinuity of $\{U_n^{\gamma};
n,\gamma\}$ is similar, and
therefore omitted.
\end{pf*}
\begin{pf*}{Proof of Lemma \protect\ref{lem02}}
We will only present the proof for the case
$\gamma= 0$. The general case follows upon minor modifications.
Denote
\[
\psi_{1,\infty}={\sup_{x_0,x\in\overline
G}\sup_j}|D\psi^{x_0}_j(x)|,\qquad
\psi_{2,\infty}={\sup_{x_0,x\in\overline G}\sup_j}|D^2\psi^{x_0}_j(x)|,
\]
and let $\varphi_{1,\infty}$, $\varphi_{2,\infty}$ be defined analogously.
Let $\varepsilon>0$ and $\bolds{\sigma}\equiv(\sigma,\xi_1,\xi
_2,\beta
,Y)\in\bolds{\Sigma}_\varepsilon$ be
given, let $\overline X=\overline X[\bolds{\sigma}]$, $\overline Y=\overline
Y[\bolds{\sigma}]$,
$\overline\tau_\rho=\overline\tau_\rho[\bolds{\sigma}]$, $\overline\tau
_G=\overline
\tau_G[\bolds{\sigma}]$,
and $\overline\tau=\overline\tau[\bolds{\sigma}]$.
Let
\[
\overline\tau_0=\inf\{t\ge\sigma\dvtx\Psi(\overline X_t)\le0\},\qquad
\overline\tau_B=\inf\{t\ge\sigma\dvtx\overline X_t\notin\mathbb{B}_{j^*}\},
\]
where we recall that $j^* = \underline{j}(\xi_1)$.
We have $\overline\tau\equiv\overline\tau_{\rho}\wedge\overline\tau_G\le\overline
\tau_{\rho}\wedge\overline\tau_0$,
because $\Psi\ge\Phi$. Also, by (\ref{34}), $\overline\tau\le\overline
\tau_B$. By It\^{o}'s
formula, for $t>\sigma$,
\begin{eqnarray*}
\Psi(\overline X_t)
&=&\Phi(\xi_2)+\int\mathbf{1}_{[\sigma,t]}(s)D\Psi(\overline
X_s)[\overline
A_s-\overline B_s]\,dW_s
\\
&&{}
+\int_\sigma^tD\Psi(\overline X_s)[\overline A_s+\overline B_s][\overline C_s+\overline D_s]\,ds
\\
&&{}
+\frac12\int_\sigma^t[\overline A_s-\overline B_s]'D^2\Psi(\overline X_s)[\overline
A_s-\overline B_s]\,ds.
\end{eqnarray*}
For $t\le\overline\tau$, using (\ref{31}) and (\ref{32}),
\begin{eqnarray*}
D\Psi(\overline X_s)[\overline A_s-\overline B_s] &=&
D\Psi(\overline X_s)[\Pi(\overline X_s)-\overline B_s]+D\Psi(\overline X_s)[\overline A_s-\Pi
(\overline X_s)]
\\
&=&-|D\Psi(\overline X_s)| (1-\alpha_s-\delta_s),
\end{eqnarray*}
where $\alpha_s=-|D\Psi(\overline X_s)|^{-1}D\Psi(\overline X_s)\cdot\overline
B_s$, and
$\delta_s=-\Pi(\overline X_s)\cdot[\overline A_s-\Pi(\overline X_s)]$. Note that
$|\alpha_s|\le1$. Moreover, using the inequality
$|\frac{v}{|v|}\cdot(\frac{u}{|u|}-\frac{v}{|v|})|\le
2|v|^{-1}|u-v|$ along with
(\ref{33}), recalling the definition of $\overline A$ and the fact $|D\Psi
(\overline X_s)|\ge1/2$, we see that
\[
|\delta_s|\le4\varepsilon\psi_{2,\infty}.
\]
Furthermore,
\begin{eqnarray*}
&&D\Psi(\overline X_s)[\overline A_s+\overline B_s][\overline C_s+\overline D_s]
\\
&&\qquad=
D\Psi(\overline X_s)[\Pi(\overline X_s)+\overline B_s][\overline C_s+\overline D_s]+D\Psi
(\overline
X_s)[\overline A_s-\Pi(\overline X_s)][\overline C_s+\overline D_s]
\\
&&\qquad=
-|D\Psi(\overline X_s)| (1+\alpha_s)(c^0+\overline D_s)+e_s,
\end{eqnarray*}
where, by (\ref{33}), for $\sigma\le t_1 \le t_2 \le\overline\tau$
\[
\int_{t_1}^{t_2}e_s\,ds\le\varepsilon+t_2-t_1.
\]
Finally, we can estimate
\[
p_s:=\tfrac12[\overline A_s-\overline B_s]'D^2\Psi(\overline X_s)[\overline A_s-\overline B_s]
\]
by $|p_s|\le2\psi_{2,\infty}$.
Shifting time by $\sigma$, we denote $\mathcaligr{G}_t=\mathcaligr
{F}_{t+\sigma}$, $\check
W_t=W_{t+\sigma}-W_\sigma$, and
\[
(\check X_t,\check D_t,\check\alpha_t,\check\delta_t,\check
e_t,\check p_t)
=(\overline X_{t+\sigma},\overline D_{t+\sigma},\alpha_{t+\sigma},\delta
_{t+\sigma}, e_{t+\sigma
},p_{t+\sigma}).
\]
Denote also $m_t=|D\Psi(\check X_s)|$, let $M$ be the
$\mathcaligr{G}_t$-martingale
\[
M_t=-\int_0^tm_s(1-\check\alpha_s-\check\delta_s)\,d\check W_s
\]
and set
%
%
\begin{eqnarray}\label{36}\quad
\mu_t:\!&=&\langle M\rangle_t=\int_0^t m_s^2(1-\check\alpha_s-\check
\delta_s)^2\,ds,
\nonumber\\
P_t&=&-\int_0^tm_s(1+\check\alpha_s)(c^0+\check D_s)\,ds,\qquad
Q_t=\int_0^t(\check e_s+\check p_s)\,ds,\\
\Psi_t &=&\Psi(\check X_t).\nonumber
\end{eqnarray}
Combining the above estimates, we have for $0\le s \le t\le\overline\tau
-\sigma$,
%
%
\begin{eqnarray}
\label{11.1}
\Psi_t
&=& \Psi_0+M_t+P_t+Q_t,
\\
\label{38}
Q_t-Q_s&\le&\varepsilon+r(t-s),
\end{eqnarray}
where $r=2\psi_{2,\infty}+1$.
Note that $m_t\ge1/2$ for $s\le\overline\tau_B-\sigma$, and recall that
$\overline\tau_B\ge\overline\tau\equiv\overline\tau_{\rho}\wedge\overline\tau
_G$. We have for $t\le
\overline\tau-\sigma$,
assuming without loss of generality $4\varepsilon\psi_{2,\infty}<1/32$,
\begin{eqnarray*}
r &=&
\frac{r}{4}  (1-\check\alpha_t-\check\delta_t+1+\check\alpha
_t+\check\delta_t)^2
\le
\frac{r}{2}(1-\check\alpha_t-\check\delta_t)^2+2r(1+\check\alpha
_t+\check\delta
_t)\\
&\le&
2rm_t^2(1-\check\alpha_t-\check\delta_t)^2+4rm_t(1+\check\alpha
_t)+2r\check\delta_t
\end{eqnarray*}
and
\[
\check\delta_t \le(1+\check\alpha_t)+\tfrac18(1-\check\alpha
_t-\check\delta_t)^2
\le2m_t(1+\check\alpha_t)+\tfrac12m_t^2(1-\check\alpha_t-\check
\delta_t)^2.
\]
Hence, for $t\le\overline\tau-\sigma$, we have
\(
r\le3rm_t^2(1-\check\alpha_t-\check\delta_t)^2+8rm_t(1+\check
\alpha_t).
\)
Thus, by~(\ref{11.1}), if $c^0$ is chosen larger than $8r$, we have
%
%
\begin{equation}
\label{15.0}
\Psi_t=\Psi_0+M_t+3r\mu_t+\widetilde P_t+\widetilde Q_t,
\end{equation}
where
\[
\widetilde P_t=-\int_0^tm_s(1+\check\alpha_s)(c^0+\check D_s-8r)\,ds,
\]
$\widetilde Q_0=0$, and
%
%
\begin{equation}\label{37}
\widetilde P_t-\widetilde P_s\le0,\qquad
\widetilde Q_t-\widetilde Q_s\le\varepsilon,\qquad  0 \le s\le t \le
\overline\tau- \sigma.
\end{equation}

We will write $\widehat{\mathbf{P}}$ for $\mathbf{P}[ \cdot
|\mathcaligr{G}_0]$, and $\widehat{\mathbf{E}}$
for the respective conditional expectation.

The proof will proceed in several steps.
\begin{Step}\label{Step1}
For some $\nu_1 \in(0, \infty)$,
\[
\sup_{\bolds{\sigma}\in\bolds{\Sigma}_{\rho_1}} \widehat{\mathbf
{E}}[(\overline\tau-\sigma
)^2] \le
\nu_1,\qquad  \mbox{a.s.}
\]
\end{Step}
\begin{Step}\label{Step2}
For some $\nu_2 \in(0, \infty)$,
\[
\sup_{\bolds{\sigma}\in\bolds{\Sigma}_{\rho_1}} \widehat{\mathbf
{E}}[(\overline\tau
_G-\sigma)^2] \le
\nu_2,\qquad  \mbox{a.s.}
\]
\end{Step}

Note that Step \ref{Step2} is immediate from Step \ref{Step1} and Lemma \ref{lem01}
because by construction, a constant control is used after time
$\overline\tau$.
\begin{Step}\label{Step3}
There exists a modulus $\vartheta_1$ such that
\[
\sup_{\bolds{\sigma}\in\bolds{\Sigma}_\varepsilon}\widehat
{\mathbf{P}}[\overline\tau_G-\sigma
>\vartheta_1(\varepsilon), \overline\tau_{\rho}>\overline\tau_G]\le
\vartheta_1(\varepsilon),\qquad \varepsilon>0.
\]
\end{Step}
\begin{Step}\label{Step4}
There exists a modulus $\vartheta_2$ such that
\[
\sup_{\bolds{\sigma}\in\bolds{\Sigma}_\varepsilon}\widehat
{\mathbf{P}}[\overline\tau_{\rho
}\le\overline\tau_G]\le\vartheta_2(\varepsilon),\qquad \varepsilon>0.
\]
\end{Step}

Based on these steps, part (i) of the lemma is established as
follows. Writing $E_\varepsilon$ for the event $\overline\tau_G-\sigma
>\vartheta_1(\varepsilon)$,
\begin{eqnarray*}
\widehat{\mathbf{E}}[\overline\tau_G-\sigma]
&=&
\widehat{\mathbf{E}}[(\overline\tau_G-\sigma)\mathbf{1}_{E_\varepsilon
}]+\widehat{\mathbf{E}}[(\overline\tau_G-\sigma)\mathbf{1}
_{E_\varepsilon^c}]\\
&\le&
 [\widehat{\mathbf{E}}[(\overline\tau_G-\sigma)^2]\widehat{\mathbf
{P}}(E_\varepsilon)
 ]^{1/2}+\vartheta_1(\varepsilon)\\
&\le&\nu_2^{1/2}[\vartheta_1(\varepsilon)+\vartheta_2(\varepsilon
)]^{1/2}+\vartheta_1(\varepsilon),
\end{eqnarray*}
where the first inequality uses Cauchy--Schwarz, and the second uses
Steps \ref{Step2}, \ref{Step3} and \ref{Step4}.

To show part (ii) of the lemma, use Steps \ref{Step3} and \ref{Step4} to write
%
%
\begin{eqnarray}\label{13.1}
\widehat{\mathbf{E}}[|\overline X-\xi_2|_{*,\overline\tau_G}^2]
&\le&\widehat{\mathbf{E}}\bigl[|\overline
X-\xi_2|_{*,\overline\tau_G}^2\mathbf{1}_{E_\varepsilon^c\cap\{\overline
\tau_{\rho}>\overline\tau_G\}}\bigr]
\nonumber\\[-8pt]\\[-8pt]
&&{}+[\vartheta_1(\varepsilon)+\vartheta_2(\varepsilon)]
\operatorname{diam}(G)^2.\nonumber
\end{eqnarray}
By (\ref{11.1}) and (\ref{38}),
we can estimate
%
%
\begin{equation}
\label{13.2}\qquad\quad
\widehat{\mathbf{E}}\Bigl[\sup_{t \in[\sigma, \overline\tau_G]}\Psi(\overline
X_{t})\mathbf{1}_{E_\varepsilon
^c\cap\{\overline\tau_{\rho}>\overline\tau_G\}}\Bigr]
\le\varphi_{1,\infty}\varepsilon+3\vartheta_1(\varepsilon
)^{1/2}\psi_{1,\varepsilon}+r\vartheta
_1(\varepsilon)+\varepsilon.
\end{equation}
Thus, noting that $\Phi(\overline X_{t})\ge0$ on $\{\overline\tau
_{\rho}>\overline\tau_G; t \in[\sigma, \overline\tau_G]\}$,
we have on this set,
\[
\Psi(\overline X_{t})\ge|\overline X_{t}-\xi_2|^2.
\]
Part (ii) of the lemma now follows on using the above inequality and
(\ref{13.2}) in (\ref{13.1}).

In order to complete the proof, we need to establish the statements
in Steps \ref{Step1}, \ref{Step3} and \ref{Step4}.
\begin{pf*}{Proof of Step \ref{Step1}}
Let $t$ be given. Let $F_t$ denote the
event $\{\Psi_s\in(0,\rho), 0\le s\le t\}$. We have
%
%
\begin{eqnarray}
\label{15.1}
\widehat{\mathbf{P}}(\overline\tau-\sigma>t)&=&\widehat{\mathbf{P}}(\overline
\tau-\sigma>t, \overline\tau_0 - \sigma> t, \overline\tau
_B-\sigma>t)\nonumber\\[-8pt]\\[-8pt]
&\le&
\widehat{\mathbf{P}}(F_t,\overline\tau_B-\sigma>t).\nonumber
\end{eqnarray}
Denote $S_u=\inf\{s\dvtx\mu_s>u\}$, where we recall that the infimum
over an empty set is taken to be $\infty$. Let $\kappa\in(0,1/16)$.
Then
%
%
\begin{eqnarray}\label{16.15}
&&\widehat{\mathbf{P}}(F_t,\overline\tau_B-\sigma>t,\mu_t>\kappa t)
\nonumber\\
&&\qquad
\le\widehat{\mathbf{P}}\bigl(\mu_t>\kappa t,
\Psi_{S_s}\in(0,\rho),  \overline\tau_B - \sigma> t,  0\le s\le
\mu_t\bigr)\nonumber\\[-8pt]\\[-8pt]
&&\qquad\le\widehat{\mathbf{P}}\bigl(\Psi_0+H_s+3rs+\widehat P_s\in(0,\rho),
0\le
s\le\kappa t\bigr)\nonumber\\
&&\qquad=\widehat{\mathbf{P}}\bigl(\widehat H_s+\widehat P_s\in(0,\rho), 0\le
s\le\kappa t\bigr),\nonumber
\end{eqnarray}
where $H$ is a standard Brownian motion (in particular, $H_s=M_{S_s}$
for $s < \mu_t$), $\widehat
P_t$~is a process that satisfies $\widehat P_s-\widehat P_u\le
\varepsilon$, $u\le
s$, and $\widehat H_s=\Psi_0+H_s+3rs$. On the event indicated in the last
line of (\ref{16.15}), one has, for every integer $k<\kappa t$, that
$\widehat H_k-\widehat H_{k-1}\ge-2$. Hence, the right-hand side of (\ref
{16.15}) can be
estimated by $m_1e^{-m_2\kappa t}$, for some positive constants
$m_1$
and $m_2$, independent of $t$ and $\kappa$, $\varepsilon$, and as a result,
%
%
\begin{equation}
\label{16.2}
\widehat{\mathbf{P}}(F_t,\overline\tau_B-\sigma>t,\mu_t>\kappa t)\le
m_1e^{-m_2\kappa t}.
\end{equation}

Next, on the event $F_t\cap\{\overline\tau_B-\sigma>t,\mu_t\le\kappa
t\}$ we have
$\int_0^tm_s^2(1-\check\alpha_s-\check\delta_s)^2\,ds\le\kappa t$, thus
\[
\int_0^t m_s^2(1-\check\alpha_s)^2\,ds\le2\kappa t
+2\int_0^tm_s^2\check\delta_s^2\,ds\le(2\kappa+16\psi_{1,\infty
}^2\psi
_{2,\infty}^2\varepsilon^2)t\le\frac t4,
\]
where we assumed without loss that
$2\kappa+16\psi_{1,\infty}^2\psi_{2,\infty}^2\varepsilon^2\le
1/4$. Consequently,
$\int_0^tm_s(1-\check\alpha_s)\,ds\le t/2$. Using $m_s\ge1/2$, we have
$\int_0^t(1+\check\alpha_s)\,ds\ge2t-t=t$, whence, letting $c^0$ be so
large that $c^0-8 r>8r$,
\[
3r\mu_t+\widetilde P_t
\le3rt-4rt=-rt,
\]
where we used $\mu_t\le\kappa t\le t$. Using (\ref{15.0}), on this
event we have $0\le\Psi_t\le\Psi_0+M_t-rt+\widetilde Q_t$. Hence, recalling
that $\widetilde Q_t\le\varepsilon$ and $\Psi_0\le\psi_{1,\infty
}\varepsilon$, denoting
$m_0=\psi_{1,\infty}+1$, and letting $\gamma$ be the ($t$-dependent)
stopping time $\gamma=\inf\{s\dvtx\mu_s>\kappa t\}$, we have, for $t\ge
r^{-1}m_0\varepsilon$,
%
%
\begin{eqnarray}
\label{17.1}
\widehat{\mathbf{P}}(F_t,\overline\tau_B-\sigma>t,\mu_t\le\kappa t)
&\le&
\widehat{\mathbf{P}}(M_t\ge rt-m_0\varepsilon, \mu_t\le\kappa t)\nonumber\\
&\le&
\widehat{\mathbf{P}}(M_{t\wedge
\gamma}\ge rt-m_0\varepsilon)\\
&\le&\frac{m_3\widehat{\mathbf{E}}[\langle M\rangle_{t\wedge\gamma
}^2]}{(rt-m
_0\varepsilon)^{4}}\le\frac{m_3
(\kappa t)^2}{(rt-m_0\varepsilon)^4}.\nonumber
\end{eqnarray}
In the second inequality above, we used the fact that $\mu_t\le\kappa
t$ implies $\gamma\ge t$, and in the third we used Burkholder's
inequality. In particular, $m_3$ does not depend on $t$ or $\kappa$
(which will allow us to use this estimate more efficiently in Step \ref{Step3}
below).
Combining (\ref{15.1}), (\ref{16.15}) and (\ref{17.1}) we obtain the
statement in Step \ref{Step1}.
\end{pf*}
\begin{pf*}{Proof of Step \protect\ref{Step3}}
We begin by observing that, from (\ref{16.15}),
%
%
\begin{equation}
\label{18.1}
\widehat{\mathbf{P}}(F_t,\overline\tau_B-\sigma>t,\mu_t>\kappa
t)\le\mathbf{P}(\widehat H_s+\widehat P_s>0, 0\le s\le\kappa t),
\end{equation}
and since $\Psi_0+\widehat P_s\le m_0\varepsilon$, this probability is
bounded by
\[
p(\varepsilon,\kappa t):=\mathbf{P}(m_0\varepsilon+H_s+3rs>0, 0\le
s\le\kappa t).
\]
The latter converges to zero as $\varepsilon\to0$ (for fixed $\kappa
$ and
$t$). Let $\overline\vartheta$ be a modulus such that
$p(\varepsilon,\overline\vartheta(\varepsilon))\le\overline\vartheta
(\varepsilon)$,
and $\frac12\overline\vartheta(\varepsilon)^{1/4}\ge m_0\varepsilon$. Taking
$t=r^{-1}(\overline\vartheta(\varepsilon))^{1/4}$ and $\kappa=r\overline
\vartheta(\varepsilon)^{3/4}$,
combining (\ref{17.1}) and (\ref{18.1}),
\[
\widehat{\mathbf{P}}\bigl(F_t,\overline\tau_B-\sigma>r^{-1}\overline\vartheta
(\varepsilon)^{1/4}\bigr)
\le\overline\vartheta(\varepsilon)+\frac{m_3\overline\vartheta(\varepsilon
)^2}{(1/2\overline
\vartheta(\varepsilon)^{1/4})^4}
=(1+16m_3)\overline\vartheta(\varepsilon).
\]
Using the above estimate in (\ref{15.1}), Step \ref{Step3} follows.
\end{pf*}
\begin{pf*}{Proof of Step \protect\ref{Step4}}
For $a>0$, let $\tau_a$ and $\tau_0$
denote the first time $[a,\infty)$, and, respectively, $(-\infty,0]$, is
hit by $\widehat H$. Since $\widehat H$ is a Brownian motion (with
drift $3r$)
starting from $\widehat H(0)\le\psi_{1,\infty}\varepsilon$, we have that
$\mathbf{P}(\tau_{\rho-\varepsilon}\le\tau_0)$ converges to zero
as $\varepsilon\to0$.
The proof is completed on noting that
\[
\widehat{\mathbf{P}}(\overline\tau_{\rho}\le\overline\tau_G)\le\mathbf
{P}(\tau_{\rho-\varepsilon}\le\tau_0),
\]
which follows from (\ref{15.0}), (\ref{37}), the relation $\widehat
H_s=\Psi_0+M_{S_s}+3r\mu_{S_s}$ for all $s < \mu_{\infty} \equiv
\sup_{t \ge0} \mu_t$ and observing that on the set where $\sigma_0 =
\sup\{S_s\dvtx s < \mu_{\infty}\} < \infty$ we have that $M_t + 3r\mu_t
= M_{\sigma_0} + 3r\mu_{\sigma_0}$, for $t \ge\sigma_0$.\qed
\noqed\end{pf*}
\noqed\end{pf*}

\section{Analysis of the game with bounded controls}\label{sec5}

The main result of this section, Theorem \ref{th2}, implies Lemma
\ref{lem3}. Fix $k,l$ such that $\min\{k, l\} \ge\max\{c^0, d^0$,
$k_{n_2},l_{n_2}\}$, where $n_2$ is as in Theorem \ref{th3plus}.
Throughout this section, $(k,l)$ will be omitted from the notation.
As in the previous section, only simple controls and strategies will be used.
Recall that
\[
\Phi(a,b,c,d;p,S)= -\tfrac12 (a-b)'S(a-b)-(c+d)(a+b)\cdot p.
\]
Fix $\gamma\in[0,1)$ and write
\begin{eqnarray*}
\Lambda^+_{\gamma}(p,S)&=&
\max_{|a|=1, 0\le c\le k}  \min_{|b|=1, 0\le d\le l}
\Phi(a,b,c,d;p,S) -\frac{\gamma^2}{2} \operatorname{Tr}(S),
\\
\Lambda^-_{\gamma}(p,S)&=&\min_{|b|=1, 0\le d\le l}  \max_{|a|=1,
0\le
c\le k}
\Phi(a,b,c,d;p,S)-\frac{\gamma^2}{2} \operatorname{Tr}(S),
\end{eqnarray*}
and consider the equations
%
%
\begin{eqnarray}\label{11}
&&\cases{\Lambda^+_{\gamma}(Du,D^2u)-h=0, &\quad in $G$,\cr
u=g, &\quad on $\partial G$,}
\\
\label{17}
&&\cases{\Lambda^-_{\gamma}(Du,D^2u)-h=0, &\quad in $G$,\cr
u=g &\quad on $\partial G$.}
\end{eqnarray}
We will write $V^{\gamma}$ and, respectively, $U^{\gamma}$ for the
functions $V^{\gamma}_{kl}$
and $U^{\gamma}_{kl}$ introduced at the beginning of Section \ref{sec4}.
\begin{theorem}\label{th2} For each $\gamma\in[0,1)$, one has the
following:

\begin{longlist}
\item
The function $U^{\gamma}$ uniquely solves (\ref{11}).
\item The function $V^{\gamma}$ uniquely solves (\ref{17}).
\end{longlist}
\end{theorem}

The proof of the theorem is based on a result on a finite time
horizon, Proposition \ref{prop1}, in which we adopt a technique of
\cite{swi}. Given a function $u\in\mathcaligr{C}(\overline G)$, $x_0\in
\overline G$,
$T\ge0$, and $Y\in M^0$, $Z\in M^0$, let
%
%
\begin{equation}\label{12}
J^{\gamma}(x_0,T,u,Y,Z)=\mathbf{E} \biggl[\int_0^{T\wedge\tau
^{\gamma
}}h(X_s^{\gamma})
\,ds+u(X^{\gamma}_{T\wedge\tau^{\gamma}}) \biggr],
\end{equation}
where $X^{\gamma}$ and $\tau^{\gamma} =\tau(x_0,Y,Z)$ are as
introduced in Section \ref{sec4}
with $X_0=x_0$, $Y=(A,C)$ and
$Z=(B,D)$.
\begin{proposition}\label{prop1}
Let $x_0\in\overline G$, $T\in[0,\infty)$ and $\gamma\in[0, 1)$.
Let $u\in
\mathcaligr{C}(\overline G)$.

\begin{longlist}
\item
If $u$ is a subsolution of (\ref{11}), then
%
%
\begin{equation}\label{13}
u(x_0)\le\sup_{\alpha\in\Gamma^0_{lk}}\inf_{Z\in
M^0_l}J^{\gamma}(x_0,T,u,\alpha[Z],Z).
\end{equation}
\item If $u$ is a supersolution of (\ref{11}), then
%
%
\begin{equation}\label{14}
u(x_0)\ge\sup_{\alpha\in\Gamma^0_{lk}}\inf_{Z\in
M^0_l}J^{\gamma}(x_0,T,u,\alpha[Z],Z).
\end{equation}
\item If $u$ is a subsolution of (\ref{17}), then
%
%
\begin{equation}\label{15}
u(x_0)\le\inf_{\beta\in\Gamma^0_{kl}}\sup_{Y\in
M^0_k}J^{\gamma}(x_0,Y,T,u,\beta[Y]).
\end{equation}
\item If $u$ is a supersolution of (\ref{17}), then
%
%
\begin{equation}\label{16}
u(x_0)\ge\inf_{\beta\in\Gamma^0_{kl}}\sup_{Y\in
M^0_k}J^{\gamma}(x_0,T,u,Y,\beta[Y]).
\end{equation}
\end{longlist}
\end{proposition}

Before proving Proposition \ref{prop1}, we show how it implies the
theorem.
\begin{pf*}{Proof of Theorem \protect\ref{th2}}
We only prove (i) since the
proof of (ii) is similar. We first argue that any solution of
(\ref{11}) must equal $U^{\gamma}$, and then show that a solution
exists. Let a solution $u$ of (\ref{11}) be given.
Fix $x_0\in\overline G$ and $\varepsilon>0$. Fix
$\alpha\in\Gamma^0_{lk}$ such that
%
%
\begin{equation}
\label{27}
U^{\gamma}(x_0)\le\inf_{Z\in M^0_l}J^{x_0}_{\gamma}(\alpha
,Z)+\varepsilon.
\end{equation}
By Proposition \ref{prop1}(ii),
%
%
\begin{equation}
\label{28}
u(x_0)\ge\inf_{Z\in M^0_l}j^{\gamma}(T,Z),
\end{equation}
where we denote
\[
j^{\gamma}(T,Z)=J^{\gamma}(x_0,T,u,\alpha[Z],Z).
\]
For the rest of the proof, we suppress $\gamma$ from the notation.
Lemma \ref{lem01} shows that there is $m_1<\infty$ such that, for
every $T\in[0,\infty)$, $\inf_{Z\in M^0_l}j(T,Z)\le m_1$.\vspace*{-2pt} Letting
$M(T)=\{Z\in M^0_l\dvtx j(T,Z)\le c_1\}$, it follows from the lower bound
on $h$ that for some $T<\infty$ that does not depend on $Z$, one has
$\mathbf{P}(\tau>T)<\varepsilon$ for all $Z\in M(T)$, where
$\tau=\tau^{x_0}(\alpha,Z)$. Fix such a $T$. Given $Z\in M(T)$, let
$\widehat Z\in M^0_l$ be equal to $Z$ on $[0,T)$, and let it assume the
constant value $(a^0,1)$ on $[T,\infty)$. Clearly $j(T,\widehat
Z)=j(T,Z)$. Also, by Lemma \ref{lem01}, denoting
$\widehat\tau=\tau^{x_0}(\alpha,\widehat Z)$, we have
\[
\mathbf{E}[(\widehat\tau-T)^+| \widehat\tau>T]\le m_2
\]
for some constant $m_2$ independent of $\varepsilon$ and $T$. By
(\ref{12}), (\ref{27}), the definition of the payoff, and using the
boundary condition $u|_{\partial G}=g$, we have for some $m_3 \in(0,
\infty)$
\begin{eqnarray*}
U(x_0)&\le& J^{x_0}(\alpha,\widehat Z)+\varepsilon\\
&\le&
J(x_0,T,u,\alpha[\widehat Z],\widehat Z)+m_3\{\mathbf{E}[(\widehat\tau
-T)^+]+\mathbf{P}
(\widehat\tau>T)\}+\varepsilon.
\end{eqnarray*}
Using $\mathbf{P}(\widehat\tau>T)=\mathbf{P}(\tau>T)$ yields
$U(x_0)\le j(T,\widehat
Z)+m_4\varepsilon=j(T,Z)+m_4\varepsilon$. Note that the infimum of
$j(T,Z)$ over $M^0_l$ is equal to that over $M(T)$. Thus, using
(\ref{28}) and sending $\varepsilon\to0$ proves that $U(x_0)\le u(x_0)$.

To obtain the reverse inequality, fix $x_0\in\overline G$. From Lemma
\ref{lem01}, there exists $m_5<\infty$ and $Z_1\in M^0_l$ such that,
for every $\alpha$,
\[
J^{x_0}(\alpha,Z_1)\le m_5.
\]
Denote $N(\alpha)=\{Z\dvtx J^{x_0}(\alpha,Z)\le m_5\}$. Clearly, for each
$\alpha$, the infimum of $J^{x_0}(\alpha,Z)$ over all $Z\in M^0_l$ is
equal to that over $Z\in N(\alpha)$. Hence,
\[
U(x_0)\ge\inf_{Z\in N(\alpha)}J^{x_0}(\alpha,Z),\qquad \alpha\in
\Gamma^0_{lk}.
\]
Using the positive lower bound on $h$ as before, it follows that
there exists a function $r\dvtx [0, \infty) \to[0, \infty)$ with $\lim
_{T\to\infty}r(T)=0$, such that
for every $\alpha$ and $Z\in N(\alpha)$ we have
$\mathbf{P}(\tau^{x_0}(\alpha,Z)>T)\le r(T)$. Therefore, for some $m_6
\in(0, \infty)$
\[
J^{x_0}(\alpha,Z)\ge J(x_0,T,u,\alpha[Z],Z)-m_6 r(T),\qquad \alpha
\in
\Gamma^0_{lk},  Z\in N(\alpha).
\]
In conjunction with Proposition \ref{prop1}(i), this shows that
$U(x_0)\ge u(x_0)- m_6 r(T)$. Since $T$ is arbitrary, we obtain
$U(x_0)\ge u(x_0)$.

Finally, we argue existence of solutions to (\ref{11}). Let
us write (\ref{11})${}_\gamma$ for (\ref{11}) with a
specific $\gamma$. For $\gamma\in(0,1)$, existence of solutions to
(\ref{11})${}_\gamma$ follows from Theorem 1.1 of \cite{CKLS}. To
handle the case $\gamma=0$, we will use the fact that any uniform
limit, as $\gamma\to0$, of solutions to (\ref{11})${}_\gamma$ is a
solution to (\ref{11})${}_0$. This fact follows by a standard
argument, that we omit. Now, since for $\gamma\in(0,1)$ we have
existence, the uniqueness statement established above shows that
$U^{\gamma}$ solves (\ref{11})${}_\gamma$. From Theorem
\ref{th3plus}, we have that the family $\{U^{\gamma}, \gamma\in
(0,1)\}$ is equicontinuous, and thus a uniform limit of solutions,
and in turn a solution to (\ref{11})${}_0$, exists.
\end{pf*}

In the rest of this section, we prove Proposition \ref{prop1}.

Let $G_n$ be a sequence of domains compactly contained in $G$ and
increasing to~$G$. Let $J_n^{\gamma}$ be defined as $J^{\gamma}$ of
(\ref{12}), with
$\tau^{\gamma}=\tau^{\gamma}(x_0,Y,Z)$ replaced by $\tau_n^{\gamma
}=\tau_n^{\gamma}(x_0,Y,Z)$, where
\[
\tau_n^{\gamma}=\inf\{t\dvtx X_t^{\gamma}\in\partial G_n\}.
\]
\begin{lemma}\label{lem6}
For every $n$, and $\gamma\in[0, 1)$ Proposition \ref{prop1} holds
with $J^{\gamma}$ replaced by
$J_n^{\gamma}$.
\end{lemma}
\begin{pf}
We follow the proof of \cite{swi}, Lemma 2.3 and Theorem 2.1.
Assume without loss that $G_0\subset\subset G_1\subset\subset
G_2\subset\subset G$. We will prove the lemma for $n=0$. Since the
claim is trivial, if $x_0\notin G_0$, assume $x_0\in G_0$. In this
proof only, write $\tau$ for $\tau_0^{\gamma}$, the exit time of
$X^{\gamma}$ from
$G_0$. Fix $\widetilde\gamma>\gamma$, let $\widetilde\tau= \tau_1^{\widetilde
\gamma
}$ and
$\sigma=\tau\wedge\widetilde\tau$. For $\varepsilon>0$, consider the
sup convolution
\[
u_\varepsilon(x)=\sup_{\xi\in{\mathbb{R}}^m} \biggl\{u(\xi)-\frac
{|\xi-x|^2}{2\varepsilon
} \biggr\},\qquad
x\in G_2,
\]
where, in the above equation only, $u$ is extended to ${\mathbb{R}}^m$ by
setting $u=0$ outside $G$. It is easy to see that there exists
$\varepsilon_0$ such that the supremum is attained inside $G$ for all
$(x,\varepsilon)\in G_2\times(0,\varepsilon_0)$. The standard mollification
$u_\varepsilon^\delta\dvtx\overline G_1\to{\mathbb{R}}$ of
$u_\varepsilon\dvtx G_2\to{\mathbb{R}}$ is well defined,
provided that $\delta$ is sufficiently small. The result \cite{swi},
Lemma 2.3, for the smooth function $u_\varepsilon^\delta$ and the
argument in
the proof of \cite{swi}, Theorem 2.1, show
\[
u_\varepsilon^\delta(x_0)\le\sup_{\alpha\in\Gamma^0}\inf_{Z\in
M^0}\mathbf{E}
\biggl[\int_0^{T\wedge\sigma}h(X^{\widetilde\gamma}_s
)\,ds+u_\varepsilon^\delta(X^{\widetilde\gamma}_{T\wedge\sigma})
\biggr]+\rho(\varepsilon,\delta
,\gamma, \widetilde\gamma),
\]
where
$\lim_{\varepsilon\to0}\lim_{\widetilde\gamma\to\gamma}\lim
_{\delta\to
0}\rho(\varepsilon,\delta,\gamma, \widetilde\gamma)=0$.
We remark here that Lemma 2.3 of \cite{swi} is written for the case
where $u$ is a subsolution of a PDE of the form
(\ref{11}) on all of $\mathbb{R}^m$ and $T\wedge\sigma$ is replaced
by $T$,
however the
proof with $u$ and $T\wedge\sigma$
as in the current setting can be carried out in exactly the same
way. Since $G_0$ is compactly contained in $G_1$, we have that for
every $\theta> 0$
\[
\sup_{\alpha\in\Gamma^0} \sup_{Z \in M^0}  \Bigl\{ \mathbf
{P} \Bigl( |T\wedge
\sigma
- T\wedge\tau| + {\sup_{0 \le s \le T}} |X^{\widetilde\gamma}_s -
X^{\gamma
}_s| >
\theta \Bigr)  \Bigr\}
\]
converges to $0$ as $\widetilde\gamma\to
\gamma$.
Moreover, $u_\varepsilon^\delta\to u_\varepsilon$ as
$\delta\to0$ and $u_\varepsilon\to u$ as $\varepsilon\to0$, where
in both cases, the
convergence is uniform on $\overline G_0$ (see ibid.). Hence, the result
follows on taking $\delta\to0$, then $\widetilde\gamma\to\gamma$ and finally
$\varepsilon\to0$.
\end{pf}
\begin{pf*}{Proof of Proposition \protect\ref{prop1}}
The main argument is similar to that of Theorem \ref{th3}, and so we
omit some of the details. We will prove only item (iv) of the
proposition, since the other items can be proved in a similar way.

Fix $x$ and $T$. Let $u$ be a supersolution of (\ref{17}).
Let $n$ be large enough so that $\operatorname{dist}(\partial G_n,
\partial G) < \rho_1$.
Write
$j_n^{\gamma}(Y,\beta)$ for $J_n^{\gamma}(x,T,u,Y,\beta[Y])$ and
$j^{\gamma}(Y,\beta)$ for
$J^{\gamma}(x,T,u,Y,\beta[Y])$. Below we will keep $\gamma$ in the
notation only if there
is scope for confusion. By Lemma \ref{lem6}, $u(x)\ge
v_n:=\inf_\beta\sup_Yj_n(Y,\beta)$, for every $n$. We need to show
$u(x)\ge v:=\inf_\beta\sup_Yj(Y,\beta)$.

Fix $\varepsilon>0$. Let $\beta_n$ be such that
%
%
\begin{equation}\label{29}
\sup_Y j_n(Y,\beta_n)\le v_n+\varepsilon
\end{equation}
and let $\tau^n_1(Y)=\tau_n^{\gamma}(x,Y,\beta_n[Y])$, $Y\in
M^0_l$. Let
$\widetilde\beta_n$ be constructed from $\beta_n$ as in the proof of
Theorem \ref{th3}, where in particular, $\beta_n[Y]$ and
$\widetilde\beta_n[Y]$ differ only on $[\tau^n_1,\infty)$, by which
$j_n(Y,\widetilde\beta_n)=j_n(Y,\beta_n)$. Choose $Y_n$ such that
\[
v\le\sup_Yj(Y,\widetilde\beta_n)\le j(Y_n,\widetilde\beta
_n)+\varepsilon
\]
and set
$\tau^n_2(Y)=\tau^{\gamma}(x,Y,\widetilde\beta_n[Y])$.
Then
\[
v-v_n-2\varepsilon\le j(Y_n,\widetilde\beta_n)-j_n(Y_n,\widetilde
\beta_n)=:\delta_n.
\]
Denote $X_n=X^x(Y_n,\widetilde\beta_n)$. Using Lemma \ref{lem02},
\[
0\le\tau^n_2-\tau^n_1<\varepsilon \quad\mbox{and}\quad
|X_n(\tau^n_1\wedge T)-X_n(\tau^n_2\wedge T)|<\varepsilon
\]
with probability tending to 1 as $n\to\infty$. It now follows from the
definition of $J_n$ and $J$ [cf. (\ref{12})] that
$\limsup_n\delta_n\le\rho(\varepsilon)$ for some modulus $\rho$. Since
$\varepsilon$ is arbitrary, this proves the result.
\end{pf*}

\section{Concluding remarks}\label{sec6}

\subsection[Identity (1.4)]{Identity (\protect\ref{44})}

Recall from (\ref{09}) that
%
%
\begin{equation}\label{48}
\Phi(a,b,c,d;p,S)= -\tfrac12 (a-b)'S(a-b)-(c+d)(a+b)\cdot p
\end{equation}
and denote
%
%
\begin{eqnarray}\label{40}
\Lambda^+(p,S)&=&
\sup_{|b|=1, 0\le d<\infty}  \inf_{|a|=1, 0\le c<\infty}
\Phi(a,b,c,d;p,S),
\\
\label{41}
\Lambda^-(p,S)&=&\inf_{|a|=1, 0\le c<\infty}  \sup_{|b|=1, 0\le
d<\infty}
\Phi(a,b,c,d;p,S)
\end{eqnarray}
[compare with (\ref{18}) and (\ref{10})].
The following proposition establishes identity (\ref{44}) that, as
discussed in the introduction, allows one to view the
infinity-Laplacian equation as a Bellman--Issacs type equation.
The result states that for the SDG of Section \ref{sec2}, the
associated Isaacs condition, $\Lambda^+=\Lambda^-$, holds. Although
we do
not make use of it in our proofs, such a condition is often invoked
in showing that the game has value (cf. \cite{FS,swi}).
\begin{proposition}\label{prop2}
For $p\in{\mathbb{R}}^m$, $p\ne0$ and $S\in\mathscr{S}(m)$,
$\Lambda^+(p,S)=\Lambda(p,S)$
and $\Lambda^-(p,S)=\Lambda(p,S)$. In particular, identity (\ref
{44}) holds.
\end{proposition}
\begin{pf}
We will only show $\Lambda^-=\Lambda$ (the proof of $\Lambda
^+=\Lambda$ being
similar). Fix $p$, $S$, and omit them from the notation.
Write $\mathcaligr{H}_k$ for $\{(a,c)\in\mathcaligr{H}\dvtx c\le k\}$ and
$\phi(y,z)$ for $\Phi(a,b,c,d)$, where $y=(a,c)$, $z=(b,d)$. Given
$\delta>0$ let
$k$ be such that
$\Lambda^-\ge\inf_{y\in\mathcaligr{H}_k}\sup_{z\in\mathcaligr{H}}\phi
(y,z)-\delta$.
Then
\[
\Lambda^-\ge\inf_{y\in\mathcaligr{H}_k}\sup_{z\in\mathcaligr
{H}_l}\phi(y,z)-\delta
=\Lambda^-_{kl}-\delta.
\]
Thus, by Lemma \ref{lem1}, $\Lambda^-\ge\Lambda$.

Next, let $\overline\phi(y) = \sup_{z \in\mathcaligr{H}} \phi(y,z)$. Fix
$\delta
\in(0, \infty)$, let $y_{\delta} = (\overline p, \delta^{-1})$, where
$\overline p=p/|p|$, and let
$z_{\delta} = (b_{\delta}, d_{\delta}) \in\mathcaligr{H}$ be such that
$\overline
\phi(y_{\delta}) \le\phi(y_{\delta}, z_{\delta}) + \delta$. Then
\begin{eqnarray*}
\Lambda^{-} &\le& \overline\phi(y_{\delta}) \le- \tfrac{1}{2}(\overline p -
b_{\delta})'S(\overline p - b_{\delta}) - (\delta^{-1} + d_{\delta})
(\overline p +
b_{\delta}) \cdot p + \delta\\
&\le&- \tfrac{1}{2}(\overline p -
b_{\delta})'S(\overline p - b_{\delta}) + \delta.
\end{eqnarray*}
Note that $b_{\delta}$ must converge to $-\overline p$ or else the middle
inequality above will say $\Lambda^{-} = - \infty$, contradicting the
bound $\Lambda^{-} \ge\Lambda$. Letting $\delta\to0$, we now have from
the third inequality that $\Lambda^{-} \le\Lambda$. The result follows.
\end{pf}

\subsection{Limit trajectory under a nearly optimal play}

In \cite{pssw}, the authors
raise questions
about the form of the limit trajectory under optimal play of the
Tug-of-War game, as the
step size approaches zero (see Section 7 therein).
It is natural to ask, similarly, whether one can characterize
(near) optimal trajectories for the SDG studied in the current
paper. Let $V$ be as given in (\ref{05}). Let $x\in\overline G$ and
$\delta>0$ be given. We say that a policy
$\beta\in\Gamma$ is $\delta$-optimal for the lower game and initial
condition $x$ if $\sup_{Y\in M}J^x(Y,\beta)\le V(x)+\delta$. When a
strategy $\beta\in\Gamma$ is given, we say that a control $Y\in M$ is
$\delta$-optimal for play against $\beta$ with initial condition $x$,
if $J^x(Y,\beta)\ge\sup_{Y'\in M}J^x(Y',\beta)-\delta$. A pair
$(Y,\beta)$ is said to be a $\delta$-optimal play for the lower game
with initial condition $x$, if $\beta$ is $\delta$-optimal for the
lower game and $Y$ is $\delta$-optimal for play against $\beta$ (both
considered with initial condition $x$).
One may ask
whether the law of the process $X^{\delta}$, under an arbitrary
$\delta$-optimal play $(\beta^{\delta},Y^{\delta})$, converges to a limit
law as $\delta\to0$; whether this limit law is the same for any choice
of such $(\beta^{\delta},Y^{\delta})$ pairs; and finally, whether an
explicit characterization of this limit law can be provided.
A somewhat less ambitious goal, that is the subject of a forthcoming
work \cite{AtBu2}
is the characterization of the limit law of $X^\delta$ under
\textit{some} choice of a $\delta$-optimal play.
The result from \cite{AtBu2} states the following.
\begin{theorem}
\label{th1new}
Suppose that $V$ is a $C^2(\overline{G})$ function and $DV \neq0$ on
$\overline{G}$. Assume there exist uniformly continuous bounded extensions,
$p$ and $q$ of $\frac{Du}{|Du|}$ and
$\frac{1}{|Du|^2}(D^2u Du-\Delta_\infty u Du)$, respectively, to
${\mathbb{R}}^m$
such that, for every $x\in{\mathbb{R}}^m$, weak uniqueness holds for
the SDE
\[
dX_t=2 p(X_t)\,dW_t+2q(X_t)\,dt,\qquad   X_0=x.
\]
Fix $x \in\overline G$ and let $X$ and $\tau$ denote such a solution
and, respectively, the corresponding exit time from $G$.
Then, given any sequence
$\{\delta_n\}_{n\ge1}$, $\delta_n \downarrow0$, there exists a
sequence of strategy-control
pairs $(\beta^{n}, Y^{n}) \in M\times\Gamma$,
$n \ge1$, with the following properties:
\begin{longlist}
\item
For every $n$, the pair $(\beta^{n},Y^{n})$ forms a
$\delta_n$-optimal play for the lower game with initial condition $x$.
\item Denoting $X^{n} = X(x, Y^{n}, \beta^{n})$ and $\tau^n
=\tau(x, Y^{n}, \beta^{n})$,
one has that $(X^{n}(\cdot\wedge\tau^n), \tau^n)$ converges in distribution
to $(X(\cdot\wedge\tau), \tau)$, as a sequence of random
variables with values in $C([0, \infty)\dvtx\overline G) \times[0,\infty]$.
\end{longlist}
An analogous result holds for the upper game.
\end{theorem}

A sufficient condition for the uniqueness to hold is that $D^2u$ is
Lipschitz on $\overline G$, since then both $ p$ and $q$ are Lipschitz,
and thus admit bounded Lipschitz extensions to ${\mathbb{R}}^m$.

%

%
\printaddresses


\begin{thebibliography}{12}

\bibitem{Aron}
%
\begin{barticle}[mr]
\bauthor{\bsnm{Aronsson},~\bfnm{Gunnar}\binits{G.}}
(\byear{1967}).
\btitle{Extension of functions satisfying {L}ipschitz conditions}.
\bjournal{Ark. Mat.}
\bvolume{6}
\bpages{551--561}.
\bid{mr={0217665}}
\end{barticle}
%
\endbibitem

\bibitem{Aron2}
%
\begin{barticle}[mr]
\bauthor{\bsnm{Aronsson},~\bfnm{Gunnar}\binits{G.}}
(\byear{1972}).
\btitle{A mathematical model in sand mechanics: {P}resentation and analysis}.
\bjournal{SIAM J. Appl. Math.}
\bvolume{22}
\bpages{437--458}.
\bid{mr={0308527}}
\end{barticle}
%
\endbibitem

\bibitem{AtBu2}
%
\begin{barticle}[vtex]
\bauthor{\bsnm{Atar},~\bfnm{R.}\binits{R.}} \AND
\bauthor{\bsnm{Budhiraja},~\bfnm{A.}\binits{A.}}
(\byear{2008}).
\btitle{On near optimal trajectories for a game associated with the
$\infty$-Laplacian.}
\bjournal{Trans. Amer. Math. Soc.}
\bvolume{360}
\bpages{77--101}.
\end{barticle}
%
\endbibitem

\bibitem{BEJ}
%
\begin{bmisc}[auto:springertagbib-v1.0]
\bauthor{\bsnm{Barron},~\bfnm{E.~N.}\binits{E.~N.}},
\bauthor{\bsnm{Evans},~\bfnm{L.~C.}\binits{L.~C.}} \AND
\bauthor{\bsnm{Jensen},~\bfnm{R.}\binits{R.}}
(\byear{2009}).
\bhowpublished{The infinity Laplacian, Aronsson's equation and their
generalizations. Preprint}.
\end{bmisc}
%
\endbibitem

\bibitem{CKLS}
%
\begin{barticle}[mr]
\bauthor{\bsnm{Crandall},~\bfnm{M.~G.}\binits{M.~G.}},
\bauthor{\bsnm{Kocan},~\bfnm{M.}\binits{M.}},
\bauthor{\bsnm{Lions},~\bfnm{P.~L.}\binits{P.~L.}} \AND
\bauthor{\bsnm{Swiech},~\bfnm{A.}\binits{A.}}
(\byear{1999}).
\btitle{Existence results for boundary problems for uniformly elliptic and
parabolic fully nonlinear equations}.
\bjournal{Electron. J. Differential Equations}
\bpages{24}.
\bid{mr={1696765}}
\end{barticle}
%
\endbibitem

\bibitem{EK}
%
\begin{bbook}[mr]
\bauthor{\bsnm{Ethier},~\bfnm{Stewart~N.}\binits{S.~N.}} \AND
\bauthor{\bsnm{Kurtz},~\bfnm{Thomas~G.}\binits{T.~G.}}
(\byear{1986}).
\btitle{Markov Processes: Characterization and Convergence}.
\bpublisher{Wiley}, \baddress{New York}.
\bid{mr={838085}}
\end{bbook}
%
\endbibitem

\bibitem{FS}
%
\begin{barticle}[mr]
\bauthor{\bsnm{Fleming},~\bfnm{W.~H.}\binits{W.~H.}} \AND
\bauthor{\bsnm{Souganidis},~\bfnm{P.~E.}\binits{P.~E.}}
(\byear{1989}).
\btitle{On the existence of value functions of two-player, zero-sum stochastic
differential games}.
\bjournal{Indiana Univ. Math. J.}
\bvolume{38}
\bpages{293--314}.
\bid{mr={997385}}
\end{barticle}
%
\endbibitem

\bibitem{Jen}
%
\begin{barticle}[mr]
\bauthor{\bsnm{Jensen},~\bfnm{Robert}\binits{R.}}
(\byear{1993}).
\btitle{Uniqueness of {L}ipschitz extensions: Minimizing the sup norm
of the
gradient}.
\bjournal{Arch. Ration. Mech. Anal.}
\bvolume{123}
\bpages{51--74}.
\bid{mr={1218686}}
\end{barticle}
%
\endbibitem

\bibitem{KoSe}
%
\begin{barticle}[mr]
\bauthor{\bsnm{Kohn},~\bfnm{Robert~V.}\binits{R.~V.}} \AND
\bauthor{\bsnm{Serfaty},~\bfnm{Sylvia}\binits{S.}}
(\byear{2006}).
\btitle{A deterministic-control-based approach to motion by curvature}.
\bjournal{Comm. Pure Appl. Math.}
\bvolume{59}
\bpages{344--407}.
\bid{mr={2200259}}
\end{barticle}
%
\endbibitem

\bibitem{pssw}
%
\begin{barticle}[mr]
\bauthor{\bsnm{Peres},~\bfnm{Yuval}\binits{Y.}},
\bauthor{\bsnm{Schramm},~\bfnm{Oded}\binits{O.}},
\bauthor{\bsnm{Sheffield},~\bfnm{Scott}\binits{S.}} \AND
\bauthor{\bsnm{Wilson},~\bfnm{David~B.}\binits{D.~B.}}
(\byear{2009}).
\btitle{Tug-of-war and the infinity {L}aplacian}.
\bjournal{J. Amer. Math. Soc.}
\bvolume{22}
\bpages{167--210}.
\bid{mr={2449057}}
\end{barticle}
%
\endbibitem

\bibitem{SoTo}
%
\begin{barticle}[mr]
\bauthor{\bsnm{Soner},~\bfnm{H.~Mete}\binits{H.~M.}} \AND
\bauthor{\bsnm{Touzi},~\bfnm{Nizar}\binits{N.}}
(\byear{2003}).
\btitle{A stochastic representation for mean curvature type geometric flows}.
\bjournal{Ann. Probab.}
\bvolume{31}
\bpages{1145--1165}.
\bid{mr={1988466}}
\end{barticle}
%
\endbibitem

\bibitem{swi}
%
\begin{barticle}[mr]
\bauthor{\bsnm{Swiech},~\bfnm{Andrzej}\binits{A.}}
(\byear{1996}).
\btitle{Another approach to the existence of value functions of stochastic
differential games}.
\bjournal{J. Math. Anal. Appl.}
\bvolume{204}
\bpages{884--897}.
\bid{mr={1422779}}
\end{barticle}
%
\endbibitem

\end{thebibliography}
\end{document}